\documentclass[12pt, a4paper]{amsart}

\usepackage{amsmath,amssymb}

\usepackage[pdftex, pdfstartview=FitH,colorlinks=true,linkcolor=blue,
urlcolor=blue,citecolor=blue,pagebackref=true]{hyperref}

\usepackage[abbrev]{amsrefs}
\usepackage{xcomment}
\usepackage{graphicx}
\usepackage{cite}
\usepackage{amsmath}
\usepackage{graphicx}
\usepackage{mathtools}

\usepackage[francais, english]{babel}

\newtheorem {theo}{Theorem}[section]
\newtheorem {lemm}[theo]{Lemma}

\newtheorem {coro}[theo]{Corollary}
\newtheorem {prop}[theo]{Proposition}

\theoremstyle{definition}

\newtheorem {defi}[theo]{Definition}

\newtheorem{rema}[theo]{Remark}

\author{Samuel Kittle} 
\email{sjmk4@cam.ac.uk}
\urladdr{www.dpmms.cam.ac.uk/person/sjmk4}
\title[Absolutely continuous self-similar measures]{Absolutely continuous self-similar measures with exponential separation}

\newcommand{\R}{\mathbb{R}}
\newcommand{\N}{\mathbb{Z}_{>0}}
\newcommand{\Z}{\mathbb{Z}}
\newcommand{\lnb}{\left\|}
\newcommand{\rnb}{\right\|}
\newcommand{\lap}{\mathop{}\!\mathbin\bigtriangleup}

\newcommand{\ceil}[1]{\left\lceil #1 \right\rceil}

\AtBeginDocument{
   \def\MR#1{}
}

\begin{document}

\begin{abstract}
In this paper, we present a sufficient condition for a self-similar measure to be absolutely continuous. In the special case of Bernoulli convolutions, we show that the Bernoulli convolution with algebraic parameter $\lambda$ is absolutely continuous provided $\lambda$ satisfies a simple condition in terms of the Mahler measure of $\lambda$, its Garsia entropy and $\lambda$. Using this, we are able to give examples of $\lambda$ for which the Bernoulli convolution with parameter $\lambda$ is absolutely continuous and for which $\lambda$ is not close to 1.
\end{abstract}

\subjclass{28A80, 60G18, 11P70}
\keywords{Bernoulli convolution, self-similar measure, absolute continuity}
\thanks{The author has received funding from the European Research Council (ERC) under the European Union’s Horizon 2020 research and innovation program (grant agreement No. 803711). The author is also supported by the CCIMI.}
\maketitle
\setcounter{tocdepth}{1}
\tableofcontents 

\section{Introduction}
\subsection{Statement of results for Bernoulli convolutions}

The main result of this paper is to give a sufficient condition for a self-similar measure to be absolutely continuous. For simplicity, we first state this result in the case of Bernoulli convolutions. First we need to define Bernoulli convolutions.
	
\begin{defi}[Bernoulli convolution] \label{defi:bernoulli_convolution}
Given some $\lambda \in (0,1)$, we define the Bernoulli convolution with parameter $\lambda$ to be the law of the random variable $Y$ given by
\begin{equation*}
Y = \sum_{n=0}^{\infty} X_n \lambda^n,
\end{equation*}
where each of the $X_n$ are i.i.d.\ random variables that have probability $\frac{1}{2}$ of being $1$ and probability $\frac{1}{2}$ of being $-1$. We denote this measure by $\mu_{\lambda}$.
\end{defi}

Bernoulli convolutions are the most well studied examples of self-similar measures which are important objects in fractal geometry. We discuss these further in Section \ref{section:lit_review}. Despite much effort, it is still not known for which $\lambda$ the measure $\mu_{\lambda}$ is absolutely continuous. The results of this paper contribute towards answering this question.
	
\begin{defi}[Mahler measure] \label{defi:mahler_measure}
Given some algebraic number $\alpha_1$ with conjugates $\alpha_2, \alpha_3, \dots, \alpha_n$ whose minimal polynomial (over $\mathbb{Z}$) has leading coefficient $C$, we define the Mahler measure of $\alpha_1$ to be
\begin{equation*}
M_{\alpha_1} = |C| \prod_{i=1}^{n}  \max \{|\alpha_i|, 1\}.
\end{equation*}
\end{defi}
	
\begin{theo} \label{theo:main_bernoulli}
Let $\lambda \in \left( \frac{1}{2}, 1 \right)$ be an algebraic number with Mahler measure $M_{\lambda}$. Suppose that $\lambda$ is not the root of any non-zero polynomial with coefficients $0, \pm 1$ and satisfies
\begin{equation}
(\log M_{\lambda} - \log 2) (\log M_{\lambda})^2 < \frac{1}{27} (\log M_{\lambda} - \log \lambda^{-1})^3 \lambda^2. \label{eq:main_bernoulli}
\end{equation}
Then the Bernoulli convolution with parameter $\lambda$ is absolutely continuous.
\end{theo}

\begin{figure}[h]
\centering
\includegraphics[scale=0.6]{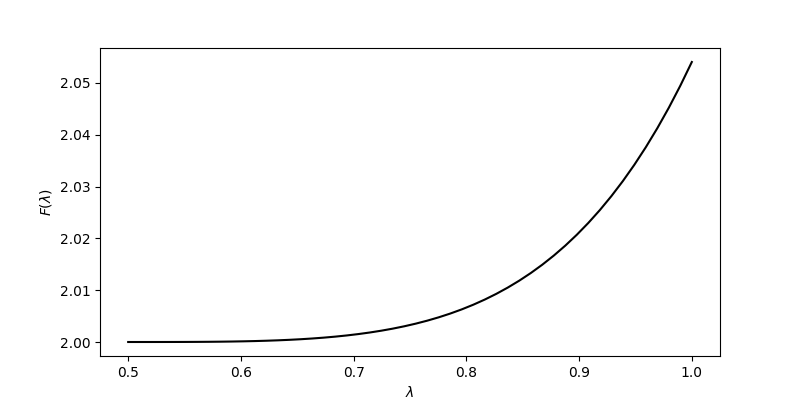}
\caption{The graph of $F$} \label{fig:fgraph}
\end{figure}

This is a corollary of a more general statement about a more general class of self-similar measures which we discuss in Section \ref{section:general_ssm}. The requirement \eqref{eq:main_bernoulli} is equivalent to $M_{\lambda} < F(\lambda)$ where $F:(\frac{1}{2}, 1) \to \R$ is some strictly increasing continuous function satisfying $F(\lambda) > 2$ and
\begin{equation*}
\left(\log F(\lambda) - \log 2 \right) \left( \log F(\lambda) \right)^2 = \frac{1}{27} \left( \log F(\lambda) - \log \lambda^{-1} \right)^3 \lambda^2
\end{equation*}
for all $\lambda \in (\frac{1}{2}, 1)$. Figure \ref{fig:fgraph} displays the graph of $F$.
	
It is worth noting that $F(\lambda) \to 2^{\frac{27}{26}} \approx 2.054$ as $\lambda \to 1$. The fact that $F(\lambda) > 2$ is important because the requirement that $\lambda$ is not the root of a polynomial with coefficients $0, \pm 1$ forces $M_{\lambda} \geq 2$ as is explained in Remark \ref{rema:mahler_bigger_than_2}.

Some parameters for Bernoulli convolutions which can be shown to be absolutely continuous using Theorem \ref{theo:main_bernoulli} are given in Table \ref{tab:examples} which can be found in Section \ref{section:examples}. The smallest value of $\lambda$ that we were able to find for which the Bernoulli convolution with parameter $\lambda$ can be shown to be absolutely continuous using this method is  $\lambda \approx  0.78207$ with minimal polynomial $X^8 - 2 X^7 - X + 1$. This is much smaller than the examples given in  \cite{VARJU_2019}, the smallest of which was $\lambda = 1 - 10^{-50 }$. We also show that for all $n \geq 8$, there is a root of the polynomial $X^n - 2X^{n-1} - X + 1$ which is in $\left( \frac{1}{2}, 1 \right)$ such that the Bernoulli convolution with this parameter is absolutely continuous.
	
\subsection{Review of existing literature} \label{section:lit_review}
For a thorough survey on Bernoulli convolutions see Peres--Schlag--Solomyak \cite{PERES_SCHLAG_WILHEM_2000} or Solomyak \cite{SOLOMYAK_2004}. For a review of recent developments see Varj\'u \cite{VARJU_2021} or Hochman \cite{MR3966837}. Bernoulli convolutions were first introduced by Jessen and Wintner in \cite{JESSEN_WINTNER_1935}. When $\lambda \in \left(0, \frac{1}{2} \right)$, it is well known that $\mu_{\lambda}$ is singular (see e.g.\ \cite{MR1507093}). When $\lambda= \frac{1}{2}$ it is clear that $\mu_{\lambda}$ is $\frac{1}{4}$ of the Lebesgue measure on $[-2,2]$. This means the interesting case is when $\lambda \in \left(\frac{1}{2}, 1 \right)$.

Bernoulli convolutions have also been studied by Erd\H{o}s. In \cite{ERDOS_1939} Erd\H{o}s showed that $\mu_{\lambda}$ is not absolutely continuous whenever $\lambda^{-1} \in (1,2)$ is a Pisot number. In his proof he exploited the property of Pisot numbers that powers of Pisot numbers approximate integers exponentially well. These are currently the only values of $\lambda \in \left( \frac{1}{2}, 1 \right)$ for which $\mu_{\lambda}$ is known not to be absolutely continuous.
	
The typical behaviour for Bernoulli convolutions with parameters in $\left( \frac{1}{2}, 1 \right)$ is absolute continuity. In \cite{ERDOS_1940} by a beautiful combinatorial argument, Erd\H{o}s showed that there is some $c < 1$ such that for almost all $\lambda \in (c,1)$, we have that $\mu_{\lambda}$ is absolutely continuous. This was extended by Solomyak in \cite{SOLOMYAK_1995} to show that we may take $c=\frac{1}{2}$. This was later extended by Shmerkin in \cite{SHMERKIN_2014} where he showed that the set of exceptional parameters has Hausdorff dimension $0$. These results have been further extended by Shmerkin in \cite{SHMERKIN_2019} who showed that for every $\lambda \in \left( \frac{1}{2}, 1 \right)$ apart from an exceptional set of zero Hausdorff dimension $\mu_{\lambda}$ is absolutely continuous with density in $L^q$ for all finite $q>1$.
	
In a ground breaking paper \cite{HOCHMAN_2014}, Hochman made progress on a related problem by showing, amongst other things, that if $\lambda \in (\frac{1}{2}, 1)$ is algebraic and not the root of any polynomial with coefficients $0, \pm 1$ then $\mu_{\lambda}$ has dimension $1$. Much of the progress in the last decade builds on the results of Hochman.
	
There are relatively few known explicit examples of $\lambda$ for which $\mu_{\lambda}$ is absolutely continuous. It can easily be shown that for example the Bernoulli convolution with parameter $2^{-\frac{1}{k}}$ is absolutely continuous when $k$ is a positive integer. This is because it may be written as the convolution of the Bernoulli convolution with parameter $\frac{1}{2}$ with another measure. Generalising this in \cite{GARSIA_1962}, Garsia showed that if $\lambda \in (\frac{1}{2}, 1)$ has Mahler measure $2$,  then $\mu_{\lambda}$ is absolutely continuous. It is worth noting that the condition that $\lambda$ has Mahler measure $2$ implies that $\lambda$ is not the root of any polynomial with coefficients $0, \pm 1$. Theorem \ref{theo:main_bernoulli} can therefore be viewed as a strengthened version of the result of Garsia.

There has also been recent progress in this area by Varj\'u in \cite{VARJU_2019}. In his paper, he showed that provided $\lambda$ is sufficiently close to $1$ depending on the Mahler measure of  $\lambda$ then $\mu_{\lambda}$ is absolutely continuous. The techniques we use in this paper are similar in many ways to those used by Varj\'u; however, we introduce several crucial new ingredients. Perhaps the most important innovation of this paper is the quantity, which we call the ``detail of a measure around a scale'', which we use in place of entropy between two scales. We discuss this further in Section \ref{section:detail}.

\subsection{Statement of results for more general self-similar measures} \label{section:general_ssm}

In this section, we discuss how the results of this paper apply to a more general class of iterated function systems.

\begin{defi}[Iterated function system] \label{defi:iterated_f}
Given some $n, d \in \N$, some homeomorphisms $S_1, S_2, \dots, S_n : \R^d \to \R^d$ and a probability vector $\mathbf{p} = (p_1, p_2, \dots, p_n)$ we say that $F = ((S_i)_{i=1}^n, \mathbf{p})$  is an iterated function system.
\end{defi}

\begin{defi}[Self-similar measure]
Given some iterated function system $F = ((S_i)_{i=1}^n, (p_i)_{i=1}^n)$ in which all of the $S_i$ are contracting similarities we say that a probability measure $\mu_F$ is a self-similar measure generated by $F$ if
\begin{equation*}
\mu_F = \sum_{i=1}^n p_i \mu_F \circ S_i ^{-1}.
\end{equation*}
\end{defi}

It is a result of J. Hutchinson in \cite[Section 3.1, Part 5]{HUTCHINSON_1981} that under these conditions there is a unique self-similar measure. Given an iterated function system, $F$, satisfying these conditions, let $\mu_F$ denote the unique self-similar measure generated by $F$. This paper only deals with a very specific class of iterated function systems.

\begin{defi}
We say that an iterated function system $F = ((S_i)_{i=1}^n, (p_i)_{i=1}^n)$ has uniform contraction ratio and uniform rotation if there is some $\lambda \in (0, 1)$, some orthogonal transformation $U: \R^d \to \R^d$ and some $a_1, a_2, \dots, a_n \in \R^d$ such that for each $i= 1, 2, \dots, n$ we have
\begin{equation*}
S_i : x \mapsto\lambda U x + a_i.
\end{equation*}
Similarly we say that the self-similar measure $\mu_F$ has uniform contraction ratio and uniform rotation when $F$ has uniform contraction ratio and uniform rotation.
\end{defi}
This notion is important because of the following lemma.
\begin{lemm} \label{lemm:sum_expression_lemma}
Let $F = ((S_i)_{i=1}^n, (p_i)_{i=1}^n)$ be an iterated function system with uniform contraction ratio and uniform rotation. Let $\lambda \in (0,1)$, let $U$ be an orthogonal transformation and let $a_1, \dots, a_n \in \R^d$ be vectors such that
\begin{equation*}
S_i : x \mapsto\lambda U x +a_i.
\end{equation*}
Let $X_0, X_1, X_2, \dots$ be i.i.d.\ random variables such that $\mathbb{P}[X_0 = a_i] = p_i$ for $i=1,\dots,n$ and let
\begin{equation*}
Y = \sum_{i=0}^{\infty} \lambda^i U^i X_i.
\end{equation*}
Then the law of $Y$ is $\mu_F$.
\end{lemm}
Using this lemma it is easy to express the self-similar measure as the convolution of many other measures. The purpose of doing this is explained in  more detail in Section \ref{section:outing_of_proof}. In order to state the main result we need the following definitions.

\begin{defi} \label{defi:v_rk}
We define the $k$-step support of an iterated function system $F$ to be given by
\begin{equation*}
V_{F,k} := \left\{ S_{j_1} \circ S_{j_2}  \circ \dots \circ S_{j_k} (0) : j_1, j_2, \dots, j_k \in \{ 1, 2, ..., n\}\right\}.
\end{equation*}
\end{defi}

\begin{defi} \label{defi:delta_fk}
Let $F$ be an iterated function system. We define the separation of $F$ after $k$ steps to be
\begin{equation*}
\Delta_{F, k} := \inf \{ |x-y| : x, y \in V_{F,k}, x \neq y \}.
\end{equation*}
\end{defi}
	
\begin{defi}
Given an iterated function system $F$ let the splitting rate of $F$, which we denote by $M_F$, be defined by
\begin{equation}
M_F := \limsup \limits_{k \to \infty} \left(\Delta_{F, k} \right) ^{-\frac{1}{k}}. \label{eq:m_f_def}
\end{equation}
\end{defi}
	
\begin{defi} \label{defi:h_f_k}
With $F$ and $X_i$ defined as above, let $h_{F,k}$ be defined by
\begin{equation*}
h_{F,k} := H\left( \sum_{i=0}^{k-1} \lambda^i U^i X_i \right).
\end{equation*}
Here $H(\cdot)$ denotes the Shannon entropy.
\end{defi}
	
\begin{defi}[Garsia Entropy] \label{defi:garsia_entropy}
Given an iterated function system $F$ with uniform contraction ratio and uniform rotation, define the \emph{Garsia entropy} of $F$ by
\begin{equation*}
h_F := \liminf_{k \to \infty} \frac{1}{k} h_{F, k}.
\end{equation*}
\end{defi}

We now have all of the definitions necessary to state the main theorem.
\begin{theo} \label{theo:main_result}
Let $F$ be an iterated function system on $\R^d$ with uniform contraction ratio and uniform rotation. Suppose that $F$ has Garsia entropy $h_F$, splitting rate $M_F$, and uniform contraction ratio $\lambda$. Suppose further that
\begin{equation*}
(d \log M_F - h_F)(\log M_F)^2 < \frac{1}{27} (\log M_F - \log \lambda^{-1})^3 \lambda^2.
\end{equation*}
Then the self-similar measure $\mu_F$ is absolutely continuous.
\end{theo}
We give examples of self-similar measures which can be shown to be absolutely continuous using this result in Section \ref{section:examples}.

\begin{rema}
Notice that it is not a requirement in the theorem for the parameters in $F$ to
be algebraic. In particular, the absolute continuity of Bernoulli convolutions
would follow even for transcendental parameters if a sufficiently good bound for the splitting rate could be proved. In Theorem \ref{theo:main_bernoulli} we bound $M_F$ for algebraic parameters using the fact that $M_{F} \leq M_{\lambda}$ which we prove in Corollary \ref{coro:garsia}. It would be interesting to bound $M_F$ for specific transcendental $\lambda$. This seems to be beyond the reach of current methods. It would also be interesting to see if the condition can be verified for almost all $\lambda \in (\frac{1}{2}, 1)$, which would allow us to recover the result of Solomyak in \cite{SOLOMYAK_1995}.
\end{rema}
\subsection{Outline of proof} \label{section:outing_of_proof}

We now describe the outline of the proof. The proof has much in common with the proof given by Varj\'u in \cite{VARJU_2019} but with some new ingredients. The most important new ingredient is the use of a new method for giving a quantitative way of measuring the smoothness of a measure at a given scale. Before defining this quantity we need to introduce the following notation.
\begin{defi}
Given an integer $d \in \N$ and some $y>0$ let $\eta_y^{(d)}$ be the density function of the multivariate normal distribution with covariance matrix $yI$ and mean $0$. Specifically let
\begin{equation*}
\eta_y^{(d)}(x) := (2 \pi y)^{-d/2} \exp \left( -\frac{1}{2y} \sum_{i=1}^d x_i^2 \right).
\end{equation*}
Where the value of $d$ is clear from context we usually just write $\eta_y$.
\end{defi}
We also use the following notation,
\begin{defi}
Given an integer $d \in \N$ and some $y>0$ let $\eta_y'$ be defined by
\begin{equation*}
\eta_y' := \frac{\partial}{\partial y} \eta_y.
\end{equation*}
This notation is only used when the value of $d$ is clear from context.
\end{defi}
We then define the following.
\begin{defi} \label{defi:detail}
Given a probability measure $\mu$ on $\R^d$ and some $r>0$ we define the \emph{detail of $\mu$ around scale $r$} by
\begin{equation*}
s_r(\mu) :=  r^2 Q(d) \lnb \mu *  \eta_{r^2}' \rnb_1
\end{equation*}
where $Q(d) := \frac{1}{2} \Gamma \left( \frac{d}{2} \right) \left( \frac{d}{2e} \right) ^{-d/2} $
\end{defi}
The factor $r^2 Q(d) $ was chosen to ensure that $s_r(\mu) \in (0, 1]$. The precise value of $Q(d)$ turns out not to matter because the factor of $Q(d)$ in Theorem \ref{theo:many_conv} ends up cancelling with the factor of $Q(d)$ in Proposition \ref{prop:entropy_to_detail}. The smaller the value of detail around a scale the smoother the measure is at that scale.

Later we show that is the detail of a measure at scale $r$ tends to $0$ sufficiently quickly as $r \to 0$ then the measure is absolutely continuous. We also show that detail decreases under convolution in a quantitative way and use this to show that the measure is absolutely continuous.

In place of $\eta_{r^2}'$, we could use another family of signed measures $\left( \nu_r \right)_{r \in \R^{+}}$ satisfying $\nu_r(\R^d) = 0$ and satisfying $\nu_{r_1}(A) = C_{r_1, r_2} \nu_{r_2}(r_2 A/ r_1)$ for every $0<r_1<r_2$ for some constant $C_{r_1, r_2}$ depending only on $r_1$ and $r_2$ for every $A \in \mathcal{B}(\R^d)$. Given such a family, we can understand something about the ``smoothness'' of $\mu$ at scale $r$ by looking at $\lnb \mu * \nu_{r} \rnb_1$. It turns out that taking $\nu_r = \eta_{r^2}'$ is a good choice because it is easy to prove Lemma \ref{lemm:suff_abs_cont} and Theorem \ref{theo:many_conv}.

First we show that provided $s_r(\mu) \to 0$ sufficiently quickly as $r \to 0$ the measure $\mu$ is absolutely continuous. Specifically we prove the following.

\begin{lemm} \label{lemm:suff_abs_cont}
Suppose that $\mu$ is a probability measure on $\R^d$ and that there exists some constant $\beta>1$ such that for all sufficiently small $r>0$ we have
\begin{equation*}
s_r(\mu) < \left( \log r^{-1} \right) ^{-\beta}.
\end{equation*}
Then $\mu$ is absolutely continuous.
\end{lemm}
This is proven in Section \ref{section:suff_abs_cont}. In order to bound the detail of the self-similar measure at a given scale we first find a quantitative bound for the detail of the convolution of many measures. Specifically we prove the following.

\begin{theo} \label{theo:many_conv}
Let $n, d \in \N$, $K>1$, $r>0$ and $\alpha_1, \alpha_2, \dots, \alpha_n \in (0,1]$. Let $m = \frac{\log n}{\log (3/2)}$. Let $\mu_1, \mu_2, \dots, \mu_n$ be probability measures on $\R^d$. Let $\alpha = \min\{ \alpha_1, \alpha_2, \dots, \alpha_n \}$.  Suppose that for all $t \in \left[ 2^{-\frac{m}{2}}r, K^m \alpha^{-m 2^m }r \right]$ and $i \in \{1, 2, \dots, n\}$ we have
\begin{equation*}
s_t(\mu_i) \leq \alpha_i.
\end{equation*}
Then we have
\begin{equation*}
s_r(\mu_1 * \mu_2 * \cdots * \mu_n) \leq C_{K, d} ^ {n-1} \alpha_1 \alpha_2 \dots \alpha_n 
\end{equation*}
where
\begin{equation*}
C_{K,d} = \frac{4}{Q(d)} \left( 1 + \frac{1}{2K^2} \right).
\end{equation*}
\end{theo}
This is be proven in Section \ref{section:decrease}. This bound is quantitatively significantly more powerful than the bound given by Varj\'u in \cite{VARJU_2019}. This is discussed further in Remark \ref{rema:like_varju}. In order to apply this theorem we need some way to express the self-similar measure as a convolution of many measures each of which have at most some detail. To do this we use entropy which is defined as follows.

\begin{defi} \label{defi:differential_entropy}
Given an absolutely continuous probability measure $\mu$ on $\R^d$ with density function $f$ we define the \emph{differential entropy} of $\mu$ by
\begin{equation*}
H(\mu) := - \int_{\R^d} f(x) \log f(x) \, dx.
\end{equation*}
Here we take $0 \log 0 = 0$.
\end{defi}

Given a continuous random variable $X$ we define $H(X)$ to be the differential entropy of its law. We also define the entropy of discrete random variables.

\begin{defi} \label{defi:discrete_entropy}
Given a discrete probability measure $\mu$ on $\R^d$ such that there is some $N \in \mathbb{Z}_{>0} \cup \{ \infty \}$,  $x_1, x_2, \dots \in \R^d$ and $p_1, p_2, \dots \in \mathbb{R}_{> 0}$ with
\begin{equation*}
\mu = \sum_{i=}^{N} p_i \delta_{x_i}
\end{equation*}
we define the \emph{Shannon entropy} of $\mu$ to be
\begin{equation}
H(\mu) := - \sum_{i=1}^{N} p_i \log p_i.
\end{equation}
\end{defi}

Similarly given a discrete random variable $X$ we define $H(X)$ to be the Shannon entropy of its law. Whenever we look at the entropy of a measure it is clear from context whether the measure is absolutely continuous or discrete. This means that using $H(\cdot)$ in both Definitions \ref{defi:differential_entropy} and \ref{defi:discrete_entropy} does not cause any problems. 

We also need the following.

\begin{defi} \label{defi:pieces}
Let $F = ((S_i)_{i=1}^n, (p_i)_{i=1}^n)$ be an iterated function system with uniform contraction ratio and uniform rotation. Suppose that $\lambda \in (0, 1)$, $U$ is an orthogonal transformation and $a_1, \dots, a_n \in \R^d$ are such that for each $i=1, 2, \dots, n$ we have
\begin{equation*}
S_i : x \mapsto \lambda U x + a_i.
\end{equation*}
Let $X_0, X_1, X_2, \dots$ be i.i.d.\ random variables such that $\mathbb{P}[X_0=a_i] = p_i$. Let $I \subset (0, \infty)$. Then we define $\mu_F^I$ to be the law of the random variable
\begin{equation*}
\sum_{i\in \Z: \lambda^i \in I} \lambda^i U^i X_i.
\end{equation*}
\end{defi}

\begin{rema}
We are only interested in the case where $I \subset (0,1]$ but allow $I \subset (0, \infty)$ to make various lemmas easier to state. We refer to the measures $\mu_F^I$ as pieces. Clearly if $I_1, I_2, \dots, I_k$ are disjoint intervals contained in $(0,1]$, then there is some measure $\nu$ such that we have
\begin{equation*}
\mu_F = \nu * \mu_F^{I_1} * \mu_F^{I_2} * \cdots * \mu_F^{I_k}.
\end{equation*}
Indeed, we can take $\nu = \mu_F^{(0,1]\backslash (I_1 \cup \dots \cup I_k)}$.
\end{rema}

We continue our outline of the proof of the main theorem. We fix a scale $r>0$ that is suitably small, but otherwise arbitrary. We aim to find suitably many disjoint intervals $I_1, I_2, \dots, I_n \subset (0, 1]$ such that $s_t (\mu_F^{I_j})$ is suitably small for $j = 1, 2, \dots, n$ for all $t$ in a suitable neighbourhood of $r$.

If we can achieve this then we can apply Theorem \ref{theo:many_conv} for the measures $\mu_F^{I_j}$ in the role of $\mu_j$. This gives us a bound on $s_r(\mu_F)$, which, if suitably good, implies the absolute continuity of $\mu_F$ via Lemma \ref{lemm:suff_abs_cont}.

In order to estimate $s_r(\mu)$ we first estimate another quantity, $H(\mu * \eta_u)$, which also measures the smoothness of the measure $\mu$. In Section \ref{section:detail_of_pieces} we prove the following result.

\begin{lemm} \label{lemm:initial_gap}
Let $F$ be an iterated function system on $\R^d$ with uniform contraction ratio and uniform rotation. Let $h_F$ be its Garsia entropy, let $M_F$ be its splitting rate, and let $\lambda$ be its contraction ratio. Then for any $M>M_F$ there is some $c>0$ such that for all $n \in \N$ we have
\begin{equation*}
H \left( \mu_F^{(\lambda^n, 1]} * \eta_{1}\right) - H\left(\mu_F^{(\lambda^n, 1]} * \eta_{M^{-2n}} \right) < (d\log M- h_F)n + c .
\end{equation*}
\end{lemm}

Under the conditions of Theorem \ref{theo:main_result}, $h_F$ is only slightly smaller than $d \log M_F$. Later we see that $H(\mu * \eta_u)$ is a non-increasing quantity in $u$. In our context this means $\left| t^2 \left. \frac{\partial}{\partial u} H(\mu_F^{(\lambda^k,1]} * \eta_u ) \right|_{u=t^2} \right|$ is small for most values of $t$ between $1$ and $M^{-n}$. Here $t^2$ is an appropriate scaling factor whose role becomes clear later.

Given a scale $s$ we can use the scaling identity
\begin{equation*}
H(\mu_F^{\lambda^kI} * \eta_{\lambda^{2k} u} ) = H(\mu_F^I * \eta_u) + dk \log \lambda
\end{equation*}
to find intervals $I$ such that $\left| s^2 \left. \frac{\partial}{\partial u} H(\mu_F^I * \eta_u) \right|_{u=s^2} \right|$ is small. We can then turn this into  an estimate for detail using the following proposition.

\begin{prop} \label{prop:entropy_to_detail}
Let $\mu$ and $\nu$ be compactly supported probability measures on $\R^d$ let $r, u$ and $v$ be positive real numbers such that $r^2=u+v$. Then
\begin{equation*}
s_r(\mu * \nu)\leq r^2 Q(d) \sqrt{ \frac{\partial }{\partial u} H(\mu * \eta_u) \frac{\partial }{\partial v} H(\nu * \eta_v) }.
\end{equation*}
\end{prop}
This is proven in Section \ref{section:entropy_det_bound}.

In Section \ref{section:main_proof}, we complete the proof of our main theorem by giving the details of the above argument to construct suitable intervals $I_j$ such that Proposition \ref{prop:entropy_to_detail}  can be applied for the measures $\mu_F^{I_j}$ and then feed the resulting estimates on detail into Theorem \ref{theo:many_conv} and finally Lemma \ref{lemm:suff_abs_cont}, as explained above. We then show that Theorem \ref{theo:main_bernoulli} follows from Theorem \ref{theo:main_result}. Finally in Section \ref{section:examples}, we give examples of self-similar measures satisfying the conditions of Theorems \ref{theo:main_bernoulli} and \ref{theo:main_result}.

\begin{rema}
Some ideas in this paper can be adapted to other settings of iterated
function systems. In a forthcoming paper \cite{KITTLE_2022}, we extend the method to study absolute continuity of Furstenberg measures.
\end{rema}

\section{Detail around a scale} \label{section:detail}
In this section we discuss the basic properties of detail around a scale. The main purpose of this section is to prove Lemma \ref{lemm:suff_abs_cont} and Theorem \ref{theo:many_conv}. 

Recall that $\eta_y' := \frac{\partial}{\partial y} \eta_y$, where $\eta_y$ is the density function of the multivariate normal distribution with mean $0$ and covariance matrix $yI$. Recall that in Definition \ref{defi:detail} we defined the detail of measure $\mu$ at scale $r$ as
\begin{equation*}
s_r(\mu) := r^2 Q(d) \lnb  \mu *  \eta_{r^2}' \rnb_1.
\end{equation*}

Detail is a quantitative measure of the smoothness of a measure at a given scale. The detail of a measure at some scale $r>0$ is close to $1$ if, for example, the measure is supported on a number of disjoint intervals of length much smaller than $r$, which are separated by a distance much greater than $r$. The detail of a measure is small if, for example, the measure is uniform on an interval of length significantly greater than $r$.

In Section \ref{section:no_increase_conv}, we prove that the detail of a probability measure does not increase if we convolve it with another probability measure. In Section \ref{section:decrease} we prove Theorem \ref{theo:many_conv}, which is a quantitative estimate on how detail decreases as we take convolutions of measures. Section \ref{section:suff_abs_cont} is devoted to the proof of Lemma \ref{lemm:suff_abs_cont} which shows that a measure is absolutely continuous provided its detail decays sufficiently fast as the scale goes to $0$. 

\begin{rema}
We motivate the definition of detail as follows. Earlier work on Bernoulli convolutions, including \cite{breuillard_varju_2020}, \cite{HOCHMAN_2014}, \cite{MR3966837}, and \cite{VARJU_2019} studied quantities like
\begin{equation*}
H(\mu * F_{r_1}) - H(\mu * F_{r_2})
\end{equation*}
where $F_r$ is a smoothing function associated to scale $r$ (for example the law of the normal distribution with standard deviation $r$ or the law of a uniform random variable on $[0,r]$). Motivated by this and the work of Shmerkin \cite{SHMERKIN_2019}, it is natural to study quantities like
\begin{equation*}
\lnb \mu * F_{r_1} \rnb_p - \lnb \mu * F_{r_2} \rnb_p.
\end{equation*}
However it turns out to be more useful to study
\begin{equation*}
\lnb \mu * ( F_{r_1} - F_{r_2}) \rnb_p
\end{equation*}
at least when $p=1$. Detail is an infinitesimal version
of this quantity with Gaussian smoothing.
\end{rema}

\subsection{No increase under convolution} \label{section:no_increase_conv}

Intuitively, convolution is a smoothing operation. This means we would not expect detail to increase under convolution. We show this in the following proposition.

\begin{prop} \label{prop:no_detail_increase}
Let $\mu$ and $\nu$ be probability measures on $\R^d$. Then we have
\begin{equation*}
s_r(\mu * \nu) \leq s_r(\mu)
\end{equation*}
\end{prop}

This is a corollary of the following Lemma.

\begin{lemm} \label{lemm:no_detail_increase}
Let $\mu$ and $\nu$ be probability measures. Then we have
\begin{equation*}
\lnb \mu * \nu *  \eta_y' \rnb_1 \leq \lnb \nu * \eta_y' \rnb_1.
\end{equation*}
Furthermore
\begin{equation}
\lnb \mu *  \eta_y' \rnb_1 \leq \lnb   \eta_y' \rnb_1 = \frac{1}{y} \cdot \frac{2}{ \Gamma \left( \frac{d}{2} \right)} \left( \frac{d}{2e} \right) ^{d/2} =\frac{1}{y Q(d)}. \label{eq:max_l1_int}
\end{equation}
\end{lemm}

\begin{rema}
It is worth noting that by \eqref{eq:max_l1_int} and the definition of detail (Definition \ref{defi:detail}) we have that $s_r(\mu) \in (0,1]$. This is the purpose of the choice of constants in Definition \ref{defi:detail}.
\end{rema}

\begin{proof}[Proof of Lemma \ref{lemm:no_detail_increase}]
For the first part simply write the measure $\nu * \eta_y' $ as $\nu * \eta_y' = \tilde{\nu}_+ - \tilde{\nu}_-$ where $\tilde{\nu}_+$ and $\tilde{\nu}_-$ are (non-negative) measures concentrated on disjoint sets. Note that this means
\begin{equation*}
\lnb \nu * \eta_y '\rnb_1 = \lnb \tilde{\nu}_+ \rnb_1 + \lnb \tilde{\nu}_- \rnb_1 
\end{equation*}
and so
\begin{align*}
\lnb \mu * \nu *  \eta_y' \rnb_1 & =  \lnb \mu * \tilde{\nu}_+ - \mu * \tilde{\nu}_-\rnb_1 \nonumber\\
& \leq  \lnb \mu \rnb_1 \lnb \tilde{\nu}_+ \rnb_1 + \lnb \mu \rnb_1 \lnb \tilde{\nu}_- \rnb_1\nonumber \\
& = \lnb \nu * \eta_y' \rnb_1.\nonumber
\end{align*}
For the second part, we need to compute
\begin{equation*}
\int_{\mathbf{x} \in \R^d} \left|  \eta_y' \right| \, d \mathbf{x} .
\end{equation*}
To do this, we work in polar coordinates. Let $s = \sqrt{ \sum_{i=1}^d x_i^2}$. Then we have
\begin{equation*}
 \eta_y'(x_1, x_2, \dots, x_d) = \left( \frac{s^2}{2y^2} - \frac{d}{2y} \right) (2 \pi y) ^ {-d/2} \exp(- \frac{s^2}{2y}) .
\end{equation*}
Noting that the $(d-1)$-dimensional surface measure of $S^{(d-1)}$ is $\frac{2 \pi^{d/2}}{\Gamma \left( d/2 \right) }$ we get
\begin{align*}
\int_{\mathbf{x} \in \R^d} \left| \eta_y' (\mathbf{x})\right| \, d \mathbf{x} &=& \frac{2 \pi^{d/2}}{\Gamma \left( \frac{d}{2} \right) } \left( -\int_{s=0}^{\sqrt{dy}} \left( \frac{s^2}{2y^2} - \frac{d}{2y} \right) (2 \pi y) ^ {-d/2} s^{d-1}\exp(- \frac{s^2}{2y}) \,ds  \right.\\ 
&&+ \left. \int_{\sqrt{dy}}^{\infty} \left( \frac{s^2}{2y^2} - \frac{d}{2y} \right) (2 \pi y) ^ {-d/2} s^{d-1}\exp(- \frac{s^2}{2y}) \,ds \right) .\nonumber
\end{align*} 
By differentiation it is easy to check that
\begin{equation*}
\int \left( \frac{s^2}{y} - d \right) s^{d-1} e^{- \frac{s^2}{2y}} \, ds = - s^d e^{- \frac{s^2}{2y}} .
\end{equation*}
Hence
\begin{align*}
\MoveEqLeft  -\int_{s=0}^{\sqrt{dy}} \left( \frac{s^2}{2y^2} - \frac{d}{2y} \right) (2 \pi y) ^ {-d/2} s^{d-1}\exp(- \frac{s^2}{2y}) \,ds\\
&\,\,\,\,\,\, +  \int_{\sqrt{dy}}^{\infty} \left( \frac{s^2}{2y^2} - \frac{d}{2y} \right) (2 \pi y) ^ {-d/2} s^{d-1}\exp(- \frac{s^2}{2y}) \,ds\\
& =   2 \cdot \frac{1}{2y} (2 \pi y) ^ {-d/2}  (dy) ^ {d/2} e^{-d/2}\\
& =  \frac{1}{y} (2 \pi)^{-d/2} d^{d/2} e^{-d/2} 
\end{align*}
which yields
\begin{equation*}
\int_{\mathbf{x} \in \R^d} \left| \eta_y' (\mathbf{x})\right| \, d \mathbf{x} = \frac{1}{y} \cdot \frac{2}{ \Gamma \left( \frac{d}{2} \right)} \left( \frac{d}{2e} \right) ^{d/2} . \qedhere
\end{equation*}
\end{proof}

\subsection{Quantitative decrease under convolution} \label{section:decrease}
In this subsection, we find a quantitative bound for the decrease of detail under convolution. Specifically we prove Theorem \ref{theo:many_conv}. We begin with a result which differs from the $n=2$ case of Theorem \ref{theo:many_conv} only in that the range of the parameter $t$ is slightly smaller.

\begin{lemm} \label{lemm:t2}
Let $\mu_1$ and $\mu_2$ be probability measures on $\R^d$, let $r>0$, $\alpha_1,\alpha_2\in (0,1]$ and let $K>1$. Suppose that for all $t \in [r/\sqrt{2}, K\alpha_1^{-\frac{1}{2}} \alpha_2^{-\frac{1}{2}}r]$ and for all $i\in \{1, 2\}$, we have
\begin{equation*}
s_t(\mu_i) \leq \alpha_i.
\end{equation*}
Then
\begin{equation*}
s_r(\mu_1*\mu_2) \leq C_{K,d} \alpha_1 \alpha_2,
\end{equation*}
where
\begin{equation*}
C_{K,d} = \frac{4}{Q(d)} \left( 1 + \frac{1}{2K^2} \right) .
\end{equation*}
\end{lemm}

We apply this lemma in the case $K \to \infty$. This means the only important property of $C_{K,d}$ is its limit as $K \to \infty$. In the case $d = 1$ this limit is $\frac{8}{\sqrt{2 e \pi}} \approx 1.93577$. We deduce Theorem \ref{theo:many_conv} from this by induction on $n$ at the end of this subsection. Before proving Lemma \ref{lemm:t2} we point out that it is analogous to \cite[Theorem 2]{VARJU_2019}.

\begin{rema} \label{rema:like_varju}
This result is similar to \cite[Theorem 2]{VARJU_2019} though more powerful. Varj\'u{}'s result states that if there is some $\alpha \in (0, \frac{1}{2})$ and some $r > 0$ such that for all $s \in \left[ \alpha^3 r, \alpha^{-3} r \right]$ we have
\begin{equation*}
1 - H(\mu ; s | 2s), 1 - H(\nu ; s | 2s) \leq \alpha
\end{equation*}
then
\begin{equation}
1 - H(\mu * \nu; r | 2r) \leq 10^8 \left(\log \alpha^{-1} \right)^3 \alpha^2. \label{eq:varju_t2}
\end{equation}
Here $H(\mu; r | 2r)$ is a quantity which Varj\'u refers to as the entropy of $\mu$ between the scales $r$ and $2r$. This quantity is always in $[0,1]$ and is closer to $1$ the smoother the measure is at scale $r$. Hence $1 - H(\mu;r|2r)$ is an analogue of $s_r(\mu)$. This result is not as powerful as Lemma \ref{lemm:t2} as it contains the factor of $\left(\log \alpha^{-1} \right)^3 $ and has a significantly larger constant term. Indeed the constant is $10^8$ instead of a constant less than $2$. Lemma \ref{lemm:t2} also has the advantage of having a significantly shorter proof and working in higher dimensions. However, note that \cite[Theorem 2]{VARJU_2019} does not follow logically from Lemma \ref{lemm:t2}.
\end{rema}

We now turn to the proof of Lemma \ref{lemm:t2}. The most important part of this proof is the following lemma.

\begin{lemm} \label{lemm:simple_int}
Let $\mu_1$ and $\mu_2$ be probability measures and let $y>0$. Then
\begin{equation*}
\lnb \mu_1 * \mu_2 * \eta_y' \rnb_1 \leq 2 \int_{\frac{y}{2}}^{\infty} \lnb \mu_1 * \eta_v' \rnb_1 \lnb \mu_2 * \eta_v' \rnb_1 \, dv .
\end{equation*}
\end{lemm}

We deduce Lemma \ref{lemm:t2} from Lemma \ref{lemm:simple_int} by simply substituting in the definition of detail. In order to prove Lemma \ref{lemm:simple_int} we need to be able to commute the $y$ derivatives. In order to do this we need the following well known result.

\begin{lemm} \label{lemm:xtoy}
Let $y>0$. Then we have
\begin{equation*}
\frac{1}{2} \lap \eta_y = \frac{\partial}{\partial y} \eta_y
\end{equation*}
where $\lap$ denotes the Laplacian
\begin{equation*}
\lap = \sum_{i=1}^d \frac{\partial^2}{\partial x_i^2} .
\end{equation*}
\end{lemm}
	
\begin{proof}
This is just a simple computation. Simply note that
\begin{equation*}
\frac{\partial}{\partial x_i} \eta_y = -\frac{x_i}{y} \eta_y
\end{equation*}
and so
\begin{equation}
\frac{1}{2} \sum_{i=1}^d \frac{\partial^2}{\partial x_i^2} \eta_y = \frac{1}{2}\left( \frac{|x|^2}{y^2} - \frac{d}{y} \right) \eta_y = \frac{\partial}{\partial y} \eta_y . \label{eq:first_euclidian_norm} \qedhere
\end{equation}
\end{proof}
In \eqref{eq:first_euclidian_norm} as in the rest of the paper we take $|\cdot|$ to be the Euclidean norm. We can now prove Lemma \ref{lemm:simple_int}. Recall the notation $\eta_y' = \frac{\partial}{\partial y} \eta_y$.

\begin{proof}[Proof of Lemma \ref{lemm:simple_int}]
First note that
\begin{equation*}
\lnb \mu * \nu * \eta_y' \rnb_1  \leq \int_y^{w} \lnb \frac{\partial}{\partial u} \left( \mu * \nu * \eta_u' \right) \rnb_1  \, du + \lnb \mu * \nu * \eta_w' \rnb_1 .
\end{equation*}
Taking $w \to \infty$ and using \eqref{eq:max_l1_int} from Lemma \ref{lemm:no_detail_increase} this gives
\begin{equation*}
\lnb \mu * \nu * \frac{\partial}{\partial y} \eta_y \rnb_1   \leq  \int_y^{\infty} \lnb \frac{\partial}{\partial u} \left( \mu * \nu * \frac{\partial}{\partial u} \eta_u\right) \rnb_1  \, du .
\end{equation*}
We can then use Lemma \ref{lemm:xtoy} and standard properties of the convolution of distributions to move the derivatives around as follows. For all $a>0$, we can write
\begin{align*}
\frac{\partial}{\partial u} \left( \mu * \nu *  \eta_u' \right) & = \frac{1}{2} \cdot \frac{\partial}{\partial u} \left( \mu * \nu * \lap \eta_u \right) \\
& =  \frac{1}{2} \cdot \frac{\partial}{\partial u} \left( \mu * \nu * \eta_{u-a} *\lap \eta_a \right) \nonumber\\
& =  \frac{1}{2} \left(\mu * \frac{\partial}{\partial u} \eta_{u-a} \right) * \left( \nu * \lap \eta_a \right) .\nonumber
\end{align*}
Letting $a= \frac{1}{2}u$ and applying Lemma \ref{lemm:xtoy} again ,this gives
\begin{equation*}
\frac{\partial}{\partial u} \left( \mu * \nu * \eta_u' \right)  =  \left(\mu * \eta_{\frac{u}{2}}' \right) * \left(\nu * \eta_{\frac{u}{2}}' \right) .
\end{equation*}
This yields
\begin{align*}
\MoveEqLeft \lnb \mu * \nu * \eta_y' \rnb_1\\
& \leq \int_{y}^{\infty} \lnb \frac{\partial}{\partial u} \left( \mu * \nu * \eta_u' \right) \rnb_1 \, du\\
& = \int_{y}^{\infty} \lnb \left(\mu * \eta_{\frac{u}{2}}' \right) * \left(\nu * \eta_{\frac{u}{2}}' \right)  \rnb_1 \, du\\
& \leq \int_{y}^{\infty} \lnb \mu * \eta_{\frac{u}{2}}' \rnb_1 \lnb \nu * \eta_{\frac{u}{2}}'  \rnb_1 \, du\\
& = 2 \int_{\frac{y}{2}}^{\infty} \lnb \mu * \eta_v' \rnb_1 \lnb \nu * \eta_v' \rnb_1 \, dv
\end{align*}
as required.
\end{proof}

We can now prove Lemma \ref{lemm:t2}.

\begin{proof}[Proof of Lemma \ref{lemm:t2}]
Using the definition of detail, applying Lemma \ref{lemm:simple_int} and using the definition of detail again we have
\begin{align}
s_r(\mu_1*\mu_2)  & =  r^2 Q(d) \lnb \mu_1*\mu_2 *  \eta_y' \rnb_1 \nonumber \\
& \leq  2 r^2 Q(d) \int_{\frac{r^2}{2}}^{\infty} \lnb \mu_1 *  \eta_v' \rnb_1 \lnb \mu_2 * \eta_v' \rnb_1 \, dv\nonumber\\
& = \frac{2 r^2}{Q(d)}  \int_{\frac{r^2}{2}}^{\infty} v^{-2} s_{\sqrt{v}}(\mu_1) s_{\sqrt{v}}(\mu_2)  \, dv. \nonumber
\end{align}
Using our assumption on detail and the fact that detail is always at most $1$, we get
\begin{align}
s_r(\mu_1*\mu_2)  & \leq  \frac{2 r^2}{Q(d)} \int_{\frac{r^2}{2}}^{K^2\alpha_1^{-1}\alpha_2^{-1}r^2} v^{-2} \alpha_1 \alpha_2 \, dv\nonumber\\
& \,\,\,\,\,\,+ \frac{2 r^2}{Q(d)} \int_{K^2\alpha_1^{-1}\alpha_2^{-1}r^2}^{\infty} v^{-2} \, dv \nonumber \\
& \leq  \frac{2 r^2}{Q(d)}  \int_{\frac{r^2}{2}}^{\infty} v^{-2} \alpha_1 \alpha_2 \, dv\nonumber\\
& \,\,\,\,\,\,  +  \frac{2 r^2}{Q(d)}  \int_{K^2\alpha_1^{-1}\alpha_2^{-1}r^2}^{\infty} v^{-2} \, dv\nonumber\\
&= \frac{2 r^2}{Q(d)} \alpha_1 \alpha_2 \left( \frac{r^{2}}{2} \right)^{-1}  +  \frac{2 r^2}{Q(d)} \left( K^2\alpha_1^{-1}\alpha_2^{-1}r^2 \right)^{-1} \nonumber \\
& =  \frac{4}{Q(d)}  \left( 1 + \frac{1}{2K^2} \right) \alpha_1 \alpha_2.\nonumber \qedhere
\end{align}
\end{proof}

We now apply Lemma \ref{lemm:t2} inductively to prove Theorem \ref{theo:many_conv}.

\begin{proof}[Proof of Theorem \ref{theo:many_conv}]
We prove this by induction. The case $n=1$ is trivial. Suppose that $n>1$. Without loss of generality we may assume that 
\begin{equation*}
0 < \alpha_1 \leq \alpha_2 \leq \dots \leq \alpha_n \leq 1
\end{equation*}
and by Lemma \ref{lemm:no_detail_increase} we may assume without loss of generality that $\alpha_i < C_{K, d}^{-1}$ for $i = 1, 2, \dots, n$. Let $n' = \left\lceil \frac{n}{2} \right\rceil$ and let $m' = \frac{\log n'}{\log (3/2)}$. Define $\nu_1, \nu_2, ..., \nu_{n'}$ and $\beta_1, \beta_2, ..., \beta_{n'}$ as follows. For $i=1, 2, ..., \left \lfloor \frac{n}{2} \right \rfloor$, let
\begin{equation*}
\nu_i = \mu_{2i-1} * \mu_{2i}
\end{equation*}
and 
\begin{equation*}
\beta_i = C_{K, d} \alpha_{2i-1} \alpha_{2i}
\end{equation*}
and if $n$ is odd, let $\nu_{n'} = \mu_n$ and $\beta_{n'} = \alpha_n$. Note that
\begin{equation}
\nu_1 * \nu_2 * \dots * \nu_{n'} = \mu_1 * \mu_2 * \dots * \mu_{n} \label{eq:measure_convolutions_the_same}
\end{equation}
and
\begin{equation}
C_{K,d} ^ {n'-1} \beta_1 \beta_2 \dots \beta_{n'} = C_{K,d}^{n-1} \alpha_1 \alpha_2 \dots \alpha_n. \label{eq:detail_target_the_same}
\end{equation}

Since $n' < n$ we just need to show that $n'$, $\left(\nu_i \right)_{i=1}^{n'}$ and $\left(\beta_i \right)_{i=1}^{n'}$ satisfy the conditions of the theorem in order to apply the inductive hypothesis. Note that $\beta_1 = \min\{\beta_1, \beta_2, \dots, \beta_{n'} \}$.
We want to use Lemma \ref{lemm:t2} to show that $s_t(\nu_i) \leq \beta_i$ for all $i = 1, 2, \dots, n'$ and for all $t \in \left[ 2^{-\frac{m'}{2}}r, K^{m'} \beta_1^{-m'2^{m'} }r \right]$. The equations \eqref{eq:measure_convolutions_the_same} and \eqref{eq:detail_target_the_same} mean that this is enough to get the required bound on $s_r(\mu_1 * \mu_2 * \dots * \mu_{n})$ by the inductive hypothesis.

To apply Lemma \ref{lemm:t2} we need to show that if $t \in \left[ 2^{-\frac{m'}{2}}r, K^{m'} \beta_1^{-m'2^{m'} }r \right]$ and $q \in \left[ 2^{-\frac{1}{2}} t, K \alpha_{2i-1} ^{-\frac{1}{2}}\alpha_{2i}^{-\frac{1}{2}} t \right]$ then $q \in  \left[ 2^{-\frac{m}{2}}r, K^m \alpha_1^{-m 2^m }r \right]$. Note that if $t \in \left[ 2^{-\frac{m'}{2}}r, K^{m'} \beta_1^{-m'2^{m'} }r \right]$ and $q \in \left[ 2^{-\frac{1}{2}} t, K \alpha_{2i-1} ^{-\frac{1}{2}}\alpha_{2i}^{-\frac{1}{2}} t \right]$ then
\begin{equation*}
q \geq 2^{-\frac{1}{2}} t \geq 2^{-\frac{m'+1}{2}} r
\end{equation*}
and
\begin{equation*}
q \leq K \alpha_{2i-1} ^{-\frac{1}{2}}\alpha_{2i}^{-\frac{1}{2}} t \leq K^{m'+1} \beta_1^{-m'2^{m'} } \alpha_1^{-1} r.
\end{equation*}
This means it is sufficient to show that
\begin{equation*}
\left[ 2^{-\frac{m'+1}{2}} r, K^{m'+1} \beta_1^{-m'2^{m'} } \alpha_1^{-1} r\right] \subset \left[ 2^{-\frac{m}{2}}r, K^m \alpha_1^{-m 2^m }r \right] .
\end{equation*}
Note that $m' +1 \leq m$ so $2^{-\frac{m'+1}{2}}r \geq 2^{-\frac{m}{2}}r$. Also we have
\begin{align*}
K^{m'+1} \beta_1^{-m'2^{m'} } \alpha_1^{-1} r & \leq  K^{m} \left( \alpha_1^{2} \right) ^{-m' 2^{m'}} \alpha_1^{-1}r]\\
& =  K^m \alpha_1^{-1 - m' 2^{m'+1}}r\\
& \leq  K^m \alpha_1^{-m 2^{m}}r
\end{align*}
as required. Hence we are done by induction.
\end{proof}
\begin{rema}
It is worth noting that the only properties of $m$ we have used are that $m\geq1$ when $n>1$ and that $m\geq m'+1$. A consequence of this is that it is possible to choose $m$ such that $m\sim \log _2 n$. It turns out that this doesn't make any difference to the bound in Theorem \ref{theo:main_result}.
\end{rema}

\subsection{Sufficiency for absolute continuity} \label{section:suff_abs_cont}
The main result of this subsection is to prove Lemma \ref{lemm:suff_abs_cont}. This lemma shows that if $s_r(\mu) \to 0$ sufficiently quickly as $r \to 0$ then $\mu$ is absolutely continuous. Lemma \ref{lemm:suff_abs_cont} follows easily from the following lemma.

\begin{lemm} \label{lemm:l1_is_abs_cont}
Let $\mu$ be a probability measure on $\R^d$ and let $y>0$. Suppose that
\begin{equation} \label{eq:norm_decay}
\int_{0^+}^y \lnb \mu *  \eta_u' \rnb_1 \, du < \infty
\end{equation}
then $\mu$ is absolutely continuous.
\end{lemm}

\begin{rema}
We use the notation $0^+$ to emphasise the fact that $\lnb \mu *  \eta_u' \rnb_1$ may not be defined at $u=0$.
\end{rema}

First we deduce Lemma \ref{lemm:suff_abs_cont} from this.

\begin{proof}[Proof of Lemma \ref{lemm:suff_abs_cont}]
Note that the requirement $s_r(\mu) < (\log r^{-1})^{-\beta}$ implies
\begin{equation*}
 r^2 Q(d) \lnb \mu *  \eta_{r^2}' \rnb_1 < (\log r^{-1})^{-\beta}.
\end{equation*}
By the conditions of Lemma \ref{lemm:suff_abs_cont} this is true for some $\beta > 1$ for all sufficiently small $r>0$. Hence there is some $y \in (0,1)$ such that we have
\begin{align}
\int_{0^+}^y \lnb \mu *  \eta_u' \rnb_1 \, du & \leq c_1 \int_{0^+}^y \frac{1}{u} \left( \log u^{-1} \right) ^{-\beta} \, du \nonumber\\
& = c_1 \int_{\log y^{-1}}^{\infty} w ^{-\beta}  \, dw\nonumber \\
& < \infty. \nonumber
\end{align}
Thus $\mu$ is absolutely continuous by Lemma \ref{lemm:l1_is_abs_cont}.
\end{proof}

We now prove Lemma \ref{lemm:l1_is_abs_cont}.
\begin{proof}[Proof of Lemma \ref{lemm:l1_is_abs_cont}]
The condition \eqref{eq:norm_decay} implies that the sequence $\mu * \eta_u$ is Cauchy as $u \to 0$ in $L^1$. This is because given some $u>v>0$, we have that
\begin{align*}
\lnb \mu * \eta_u - \mu * \eta_v \rnb_1 & \leq \int_{v}^{u} \lnb \mu * \eta_w' \rnb_1 \, dw\\
& \leq \int_{0^+}^{u} \lnb \mu * \eta_w' \rnb_1 \, dw\\
& \to 0.
\end{align*}
Since the space $L^1$ is complete, there is some absolutely continuous measure $\tilde{\mu}$ such that $\mu * \eta_u \to \tilde{\mu}$ with respect to $L^1$ as $u \to 0$. We now just need to check that $\mu = \tilde{\mu}$.
	
Suppose for contradiction that $\tilde{\mu} \neq \mu$. The set of open subsets of $\R^d$ is a $\pi$-system generating $\mathcal{B}(\R^d)$. Therefore there is some open set $U \subset \R^d$ such that 
\begin{equation*}
\mu(U) \neq \tilde{\mu}(U) .
\end{equation*}
		
We assume for simplicity that
\begin{equation*}
\mu(U) > \tilde{\mu}(U) .
\end{equation*}
The opposite case is almost identical and we leave it to the reader. By regularity, there exists some compact set $K \subset U$ such that
\begin{equation*}
\mu(K) > \tilde{\mu}(U) .
\end{equation*}
Let $\varepsilon = \min \{ \text{dist}(K, U^C), \mu(K) -  \tilde{\mu}(U)\}$. We now consider $\mu * (\eta_u|_{B_{\varepsilon}})$ where $B_{\varepsilon}$ is the ball of radius $\varepsilon$ centred at $0$. We have
\begin{align*}
(\mu * \eta_u)(U) & \geq (\mu * (\eta_u|_{B_{\varepsilon}}))(U)\\
& \geq \lnb \eta_u|_{B_{\varepsilon}} \rnb_1 \mu(K)\\
& \geq \lnb \eta_u|_{B_{\varepsilon}} \rnb_1 (\tilde{\mu}(U) + \varepsilon)\\
& \to \tilde{\mu}(U) + \varepsilon
\end{align*}
as $u \to 0$. This contradicts the requirement
\begin{equation*}
\left(\mu* \eta_u \right) (U) \to \tilde{\mu}(U) 
\end{equation*}
as $u \to 0$. This shows that $\mu = \tilde{\mu}$ and so, in particular, $\mu$ is absolutely continuous.
\end{proof}

\section{Bounding detail using entropy} \label{section:entropy_det_bound}
The purpose of this section is to prove Proposition \ref{prop:entropy_to_detail}, which estimates the detail of a convolution of measures in terms of the quantity $\frac{\partial}{\partial u} H(\mu * \eta_u)$ for both convolution factors in the role of $\mu$.

The most important ingredient in proving Proposition \ref{prop:entropy_to_detail} is the following proposition.

\begin{prop} \label{prop:tmpc1}
Let $\mu$ be a probability measure on $\R^d$ with finite variance and let $y>0$. Then we have
\begin{equation*}
\frac{1}{2} \lnb \nabla \mu * \eta_y \rnb_1^2 \leq \frac{\partial}{\partial y} H(\mu * \eta_y).
\end{equation*}
\end{prop}

This proposition is the reason for the estimate in Proposition \ref{prop:entropy_to_detail} to be an estimate on the detail of a convolution of two measures rather than an estimate on the detail of one measure. This is because we use Lemma \ref{lemm:xtoy} to estimate $\lnb \mu * \nu * \eta_y' \rnb_1$ in terms of $\lnb \nabla \mu * \eta_u \rnb_1$ and $\lnb \nabla \nu * \eta_v \rnb_1$. 

To prove this proposition we use Fisher information.

\begin{defi}[Fisher information]
Let $\mu$ be an absolutely continuous probability measure on $\R^d$. Let $f$ be the density function of $\mu$. Suppose that $f$ is smooth. Then we define the \emph{Fisher information} of $\mu$ by
\begin{equation*}
J(\mu) := \int_{\R^d}  \frac{\left| \nabla f (x) \right|^2}{f(x)}   \, dx.
\end{equation*}
\end{defi}

\begin{theo}[de Bruijn's identity] \label{theo:de_bruijn}
Let $\mu$ be a probability measure on $\R^d$ with finite variance and let $y>0$. Then we have
\begin{equation*}
\frac{\partial}{\partial y} H(\mu * \eta_y) = \frac{1}{2} J(\mu * \eta_y).
\end{equation*}
In particular, the derivative on the left exists for all $y > 0$.
\end{theo}
\begin{proof}
This is proven in for example \cite[Theorem C.1]{Johnson_2004}.
\end{proof}

\begin{proof}[Proof of Proposition \ref{prop:tmpc1}]
Let $f$ be the density function of $\mu * \eta_y$. Note that we define
\begin{equation*}
\lnb \nabla f \rnb_1 := \int_{\R^d} |\nabla f(x)| \, dx
\end{equation*}
where $| \cdot |$ denotes the Euclidean norm. Note that we have
\begin{align*}
\lnb \nabla f \rnb_1 = \int_{\R^d} |\nabla f(x)| \, dx 
= \int_{\R^d}  \frac{ |\nabla f(x)|}{f(x)}  f(x)\, dx
\end{align*}
and so by Jensen's inequality
\begin{align*}
\lnb \nabla f \rnb_1^2  = \left(\int_{\R^d}  \frac{ |\nabla f(x)|}{f(x)}  f(x)\, dx \right)^2
 \leq \int_{\R^d} \left( \frac{ |\nabla f(x)|}{f(x)} \right)^2 f(x)\, dx
 = J(\mu * \eta_y).
\end{align*}
The result now follows by Theorem \ref{theo:de_bruijn}.
\end{proof}

We are now ready to prove Proposition \ref{prop:entropy_to_detail}.
\begin{proof}[Proof of Proposition \ref{prop:entropy_to_detail}]
Let $y = r^2$ and let $u, v>0$ be such that $u+v = r^2$. First note that by Lemma \ref{lemm:xtoy}, we have
\begin{align*}
\mu * \nu * \eta_y'(x) & = \frac{1}{2} \sum_{i=1}^d \frac{\partial^2}{\partial x_i^2} \mu * \nu * \eta_y(x) \\
& = \frac{1}{2} \sum_{i=1}^d \frac{\partial}{\partial x_i} \mu* \eta_u * \frac{\partial}{\partial x_i} \nu * \eta_v (x)\\
& = \frac{1}{2} \int_{\R^d} \sum_{i=1}^d \left( \frac{\partial}{\partial x_i} \mu * \eta_u (x-a) \right) \left( \frac{\partial}{\partial x_i} \nu * \eta_v (a) \right) \, da.
\end{align*}
In particular, by Cauchy-Schwartz
\begin{align*}
\left| \mu * \nu * \eta_y'(x) \right| \leq \frac{1}{2} \int_{\R^d} \left| \nabla \mu * \eta_u(x-a) \right| \cdot \left| \nabla \nu * \eta_v(a) \right| \, da
\end{align*}
and so
\begin{align*}
\lnb \mu * \nu * \eta_y' \rnb_1 \leq \frac{1}{2} \lnb \nabla \mu * \eta_u \rnb_1 \cdot \lnb \nabla \nu * \eta_v \rnb_1.
\end{align*}
By Proposition \ref{prop:tmpc1}, we then have
\begin{align*}
\lnb \mu * \nu * \eta_y' \rnb_1 \leq \sqrt{\frac{\partial}{\partial u} H(\mu * \eta_u) \frac{\partial}{\partial v} H(\nu * \eta_v)}
\end{align*}
and so by the definition of detail
\begin{equation*}
s_r(\mu * \nu) \leq r^2 Q(d) \sqrt{\frac{\partial}{\partial u} H(\mu * \eta_u) \frac{\partial}{\partial v} H(\nu * \eta_v)},
\end{equation*}
as required.
\end{proof}

\section{Entropy of pieces} \label{section:detail_of_pieces}

The purpose of this Section is to prove Lemma \ref{lemm:initial_gap} which provides an estimate for the difference of the entropy of $\mu_F^{(\lambda^k, 1]}$ smoothed at two appropriate scales in terms of the Garsia entropy of the iterated function system $F$. We now recall the definition of $\mu_F^I$ from Definition \ref{defi:pieces}. Let $F = \left( (S_i)_{i=1}^n, (p_i)_{p=1}^n \right)$ be an iterated function system such that there is some orthogonal $U$ and some $\lambda \in (0, 1)$ and $a_1, a_2, \dots, a_n \in \R^d$ such that
\begin{equation*}
S_i: x \mapsto \lambda U x + a_i.
\end{equation*}
Let $I \subset (0, \infty)$. Then we define $\mu_F^{I}$ to be the law of the random variable
\begin{equation*}
\sum_{n \in \Z : \lambda^n \in I} \lambda^n U^n X_i
\end{equation*}
where the $X_i$ are i.i.d.\ random variables with $\mathbb{P}[X_i = a_i] = p_i$. The purpose of this subsection is to prove the following.

\begin{lemm} \label{lemm:entropy_small_overlap_lower_bound}
Let $n \in \N$, $r,R \in \R_{>0}$. Let $x_1, \dots, x_n \in \R^d$ be such that $|x_i-x_j| \geq 2R$ for $i \neq j$. Let $\mathbf{p} = (p_1, p_2 , \dots, p_n)$ be a probability vector and let
\begin{equation*}
\mu =  \sum_{i = 1}^n p_i \delta_{x_i}.
\end{equation*}
Then
\begin{equation*}
H(\mu * \eta_{r^2} ) \geq 	d\log r +  H (\mathbf{p}) -c
\end{equation*}
for some constant $c$ depending only on $d$ and the ratio $R/r$.
\end{lemm}

Here and throughout the paper $H(\mathbf{p})$ means $-\sum_{i=1}^n p_i \log p_i$ and in the case where $\mathbf{p}$ has infinitely many components we take  $H(\mathbf{p})$ to be $-\sum_{i=1}^{\infty} p_i \log p_i$. This lemma is unsurprising. This is because if we had some other measure $\nu$ supported on a ball of radius $R$ centred at $0$ then $H(\mu * \nu) = H(\nu) + H(\mathbf{p})$. The overlaps of some parts of the normal distributions means that $H(\mu * \eta_{r^2})$ is slightly less than $H(\eta_{r^2}) + H(\mathbf{p})$. We show that this difference is only some constant. This is sufficient as $H(\eta_{r^2}) = d \log r + c$. We will leave the proof of Lemma \ref{lemm:entropy_small_overlap_lower_bound} until later in the section.

\begin{lemm} \label{lemm:garsia_dec}
Let $k \in \N$. Then $H(\mu_F^{(\lambda^k, 1]}) \geq k h_F$.
\end{lemm}

\begin{proof}[Proof of Lemma \ref{lemm:garsia_dec}]
Note that $H(\mu_F^{(\lambda^n, 1]}) = h_{F, n}$ with $h_{F, n}$ as in Definition \ref{defi:h_f_k} and $h_F := \liminf_{k \to \infty} \frac{1}{k} h_{F, k}$ and that $h_{F, k} := H\left(\sum_{i=0}^{k-1} \lambda^i U^i X_i \right)$. Note that we have $h_{F, a+b} \leq h_{F,a} + h_{F, b}$. This is because $\sum_{i=0}^{a+b-1} \lambda^i U^i X_i$ is a function of $\sum_{i=0}^{a-1} \lambda^i U^i X_i$ and $\sum_{i=a}^{a+b-1} \lambda^i U^i X_i$ and

\begin{equation*}
H\left(\sum_{i=a}^{a+b-1} \lambda^i U^i X_i\right) =  H\left(\sum_{i=0}^{b-1} \lambda^i U^i X_i\right).
\end{equation*}

Suppose for contradiction there is some $k$ such that $h_{F, k} < k h_F$. Then we have $\frac{1}{ak} h_{F, ak} \leq \frac{1}{k} h_{F, k} < h_F$ for all $a \in \N$. This contradicts the definition of $h_F$.
\end{proof}

\begin{lemm} \label{lemm:ent_increase}
Suppose that $X$ and $Y$ are random variables with finite entropy either both discrete or both absolutely continuous. Then
\begin{equation*}
H(X+Y) \geq H(X)
\end{equation*}
\end{lemm}
\begin{proof}
This is well known. See for example \cite[Lemma 1.15]{Johnson_2004}.
\end{proof}
\begin{coro} \label{coro:ent_monotone}
Suppose that $I_1 \subset I_2$. Then
\begin{equation*}
H(\mu_F^{I_1}) \leq H(\mu_F^{I_2})
\end{equation*}
\end{coro}
\begin{proof}
This follows immediately from Lemma \ref{lemm:ent_increase} and the definition of $\mu_F^I$.
\end{proof}
This is sufficient to prove Lemma \ref{lemm:initial_gap} as shown below.

\begin{proof}[Proof of Lemma \ref{lemm:initial_gap}]
Note that provided $n$ is sufficiently large we have $\Delta_{F, n} \geq  M^{-n}$. In other words $\mu_F^{(\lambda^n, 1]}$ is supported on a number of points each of which are separated by a distance of at least $M^{-n}$. By Lemma \ref{lemm:garsia_dec} we also have  that $H(\mu_F^{(\lambda^n, 1]}) \geq n h_F$. Hence by Lemma \ref{lemm:entropy_small_overlap_lower_bound}, we have
\begin{equation*}
H(\mu_F^{(\lambda^n, 1]} * \eta_{M^{-2n}}) \geq n h_F -d n \log M - c.
\end{equation*}
We also have by Corollary \ref{coro:ent_monotone} that $H(\mu_F^{(\lambda^n, 1]} * \eta_{1}) < H(\mu_F * \eta_{1}) < \infty$. This gives the required result.
\end{proof}

To prove Lemma \ref{lemm:entropy_small_overlap_lower_bound}, we need to introduce the following.

\begin{defi} \label{defi:nonprob_ent}
Given a finite measure $\mu$ it is convenient to define
\begin{equation*}
H(\mu) := \lnb \mu \rnb_1 H\left(\frac{\mu}{\lnb \mu \rnb_1}\right) .
\end{equation*}
\end{defi}

\begin{lemm} \label{lemm:concave}
Let $\mu_1, \mu_2, \dots$ be absolutely continuous finite measures on $\R^d$ with finite differential entropy such that $\sum_{i=1}^{\infty} \lnb \mu_i \rnb_1 < \infty$ and both $H \left( \sum_{i=N}^{\infty} \mu_i \right)$ and $\sum_{i=N}^{\infty} H \left(  \mu_i \right)$ tend to $0$ as $N \to \infty$. Then we have
\begin{equation}
H \left( \sum_{i=1}^{\infty} \mu_i \right) \geq \sum_{i=1}^{\infty} H \left(  \mu_i \right). \label{eq:concave_entropy_lemma}
\end{equation}
\end{lemm}

\begin{proof}
First we wish to show that if $\mu$ and $\nu$ are finite measures with finite entropy then

\begin{equation}
H(\mu + \nu) \geq H(\mu) + H(\nu). \label{eq:concave_thing}
\end{equation}
Define the function $h$ by
\begin{align*}
h : [0,\infty) & \to \R\\
x & \mapsto - x \log x.
\end{align*}
Note that $h$ is concave. Let $\mu$ and $\nu$ have density functions $f$ and $g$ respectively. Note that we have

\begin{align*}
H(\mu + \nu) & = (\lnb \mu \rnb_1 + \lnb \nu \rnb_1) \int_{\R^d} h \left( \frac{f+g}{\lnb \mu \rnb_1 + \lnb \nu \rnb_1} \right)\\
& \geq (\lnb \mu \rnb_1 + \lnb \nu \rnb_1) \int_{\R^d} \frac{\lnb \mu \rnb_1}{\lnb \mu \rnb_1 + \lnb \nu \rnb_1} h\left( \frac{f(x)}{\lnb \mu \rnb_1} \right) + \frac{\lnb \nu \rnb_1}{\lnb \mu \rnb_1 + \lnb \nu \rnb_1} h\left( \frac{g(x)}{\lnb \mu \rnb_1} \right) \, dx\\
& = H(\mu) + H(\nu)
\end{align*}

as required. Applying \eqref{eq:concave_thing} inductively gives

\begin{equation*}
H \left( \sum_{i=1}^{N} \mu_i \right) \geq \sum_{i=1}^{N} H \left(  \mu_i \right).
\end{equation*}

Putting $\sum_{i=N}^{\infty} \mu_i$ in the role of $\mu_N$ and noting that $H \left( \sum_{i=N}^{\infty} \mu_i \right)$ and $\sum_{i=N}^{\infty} H \left(  \mu_i \right)$ tend to $0$ as $N \to \infty$ gives \eqref{eq:concave_entropy_lemma} as required.

\end{proof}

\begin{lemm} \label{lemm:max_inc}
Let  $\mathbf{p} = (p_1, p_2, \dots)$ be a probability vector and let $\mu_1, \mu_2, \dots$ be either all be absolutely continuous measures with finite entropy or all be discrete measures such that $\lnb \mu_i \rnb_1 = p_i$. Then
\begin{equation}
H\left(\sum_{i=1}^{\infty} \mu_i\right) \leq H(\mathbf{p}) + \sum_{i=1}^{\infty} H(\mu_i). \label{eq:inf_h_entropy}
\end{equation}
In particular if $p_i = 0$ for all $i > k$ for some $k\in \N$, then
\begin{equation}
H\left(\sum_{i=1}^k \mu_i\right) \leq \sum_{i=1}^k H(\mu_i) + \log k . \label{eq:k_h_entropy}
\end{equation}
\end{lemm}

\begin{proof}
First we prove \eqref{eq:inf_h_entropy}. To begin with we deal with the case that the measures are all absolutely continuous. Let the density function of $\mu_i$ be $f_i$. Let $h:[0, \infty) \to \R$ be defined as in the proof of Lemma \ref{lemm:concave}. Using the fact that $\sum_{i=1}^{\infty} \mu_i$ is a probability measure and the inequality $h\left( \sum_{i=1}^{\infty} a_i \right) \leq \sum_{i=1}^{\infty} h(a_i) $ we get
\begin{align}
H\left(\sum_{i=1}^{\infty} \mu_i\right) & =  \int_{\R^d} h\left(\sum_{i=1}^{\infty} f_i\right) \label{eq:max_inc1} \\
& \leq  \sum_{i=1}^{\infty}\int_{\R^d} h(f_i) \label{eq:max_inc2}\\
& = \sum_{i=1}^{\infty} \int_{\R^d} \left( - f_i(x) \log(p_i^{-1} f_i) - f_i(x) \log p_i \right)  \, dx\label{eq:max_inc3} \\
& = \sum_{i=1}^{\infty} \int_{\R^d} p_i H(p_i^{-1} f_i(x)) \, dx + h(p_i)   \label{eq:max_inc4} \\
&= \sum_{i=1}^{\infty} H(\mu_i) + H(\mathbf{p}) \nonumber.
\end{align}
The case where the measures are discrete follows from taking the density function of the measures and the integrals in \eqref{eq:max_inc1}, \eqref{eq:max_inc2}, \eqref{eq:max_inc3} and \eqref{eq:max_inc4} to be with respect to the counting measure on $\R^d$ rather than with respect to the Lebesgue measure.

For \eqref{eq:k_h_entropy} we simply apply \eqref{eq:inf_h_entropy} with $p_i = 0$ for $i > k$. We note that this gives $H(\mathbf{p}) \leq \log k$.
\end{proof}
	
\begin{lemm} \label{lemm:sep_inc}
Let $\mu$ and $\nu$ be probability measures on $\R^d$. Suppose that $\mu$ is a discrete measure supported on finitely many points with separation at least $2R$ and that $\nu$ is an absolutely continuous measure with finite entropy whose support is contained in a ball of radius $R$.  Then
\begin{equation*}
H(\mu * \nu) = H(\mu) + H(\nu) .
\end{equation*}
\end{lemm}

\begin{proof}
Let $n \in \N$, $p_1, p_2, \dots, p_n \in (0, 1)$ and $x_1, x_2, \dots, x_n \in \R^d$ be chosen such that
\begin{equation*}
\mu = \sum_{i = 1}^n p_i \delta_{x_i}.
\end{equation*}

Let $f$ be the density function of $\nu$. Note that the density function of $\mu * \nu$, which we denote by $g$,  can be expressed as

\begin{equation*}
g(x)  = 
\begin{cases}
p_i f(x - x_i)  & |x_i - x| < R \text{ for some } i,\\
0 & \text{otherwise}.\\
\end{cases}.
\end{equation*}

We then compute

\begin{align*}
H(\mu * \nu) & =  \sum_{i=1}^{n} \int_{B_R(x_i)} - g(x) \log \, g(x) \, dx\\
& = \sum_{i=1}^{n} \int_{B_R(0)} - p_i f(x) \log \left( p_i f(x) \right) \, dx \\
& =  \sum_{i=1}^{n} \int_{B_R(0)} -  p_i f(x) \log \left( f(x) \right) \, dx \\
& \,\,\,\,\,\, +  \sum_{i=1}^{n} \int_{B_R(0)} -  p_i f(x) \log \left( p_i\right) \, dx \\
& =  H(\mu) + H(\nu) \qedhere
\end{align*}
\end{proof}

We are now ready to prove Lemma \ref{lemm:entropy_small_overlap_lower_bound}.

\begin{proof}[Proof of Lemma \ref{lemm:entropy_small_overlap_lower_bound}]
Given $k \in \Z_{\geq 2}$ define
\begin{equation*}
\tilde{\eta}_k := \eta_{r^2} | _ {A_{\frac{(k-2)R}{\sqrt{d}}, \frac{(k-1)R}{\sqrt{d}}}}
\end{equation*}
where $A_{a,b} := \{x \in \R^d : |x| \in [a, b)\}$. 

We now wish to write $\mu$ as the sum of $k^d$ measures each of which are supported on points separated by at least $\frac{2(k-1)R}{\sqrt{d}}$. Given $\mathbf{m} \in \Z^d$, define
\begin{equation*}
B_{\mathbf{m}} := \left\{ x \in \R^d : x \in \mathbf{m} + [0, 1)^d \right\},
\end{equation*}
and given $\mathbf{j} \in \left( \Z / k\Z \right)^d$ we define
\begin{equation*}
\tilde{B}_{\mathbf{j}} := \bigcup_{\mathbf{m} \in \Z^d : \mathbf{m} \equiv \mathbf{j}} B_{\mathbf{m}}.
\end{equation*}

Now given $k \in \Z_{\geq 2}$ and $\mathbf{j} \in \left( \Z / k\Z \right)^d$ we define
\begin{equation*}
\nu_{\mathbf{j}, k} := \sum_{i: x_i \in \frac{2R}{\sqrt{d}} \tilde{B}_{\mathbf{j}}} p_i \delta_{x_i}.
\end{equation*}
Note that given any $k \in \Z_{\geq 2}$ we have
\begin{equation*}
\mu = \sum_{\mathbf{j} \in \left( \Z / k\Z \right)^d} \nu_{\mathbf{j}, k}.
\end{equation*}

Note that if $x_i$ and $x_j$ are distinct points in the support of $\mu$ then there cannot be any $\mathbf{m} \in \Z^d$ such that $x_i, x_j \in \frac{2R}{\sqrt{d}} B_{\mathbf{m}}$ as this would contradict the requirement $|x_i-x_j| > 2R$. If in addition, $x_i$ and $x_j$ are in the support of $\nu_{\mathbf{j}, k}$ for some $\mathbf{j} \in \Z^d$ and $k \in \Z_{\geq 2}$ then the distance between $x_i$ and $x_j$ must be at least $\frac{2(k-1)R}{\sqrt{d}}$.

By Lemma \ref{lemm:max_inc} we have 
\begin{align*}
\sum_{\mathbf{j} \in \left( \Z / k\Z \right)^d}  H(\nu_{\mathbf{j}, k}) &\geq H(\mu) - d \log k \\ 
&=  H(\mathbf{p}) - d \log k .
\end{align*}
Also by Lemma \ref{lemm:sep_inc}
\begin{align*}
H(\nu_{\mathbf{j},k} * \tilde{\eta}_k) & = \lnb \nu_{\mathbf{j},k} \rnb_1  \lnb \tilde{\eta}_k \rnb_1 H \left( \frac{\nu_{\mathbf{j},k} * \tilde{\eta}_k}{\lnb \nu_{\mathbf{j},k} \rnb_1 \lnb \tilde{\eta}_k \rnb_1 } \right)\\
& = \lnb \nu_{\mathbf{j},k} \rnb_1 \lnb \tilde{\eta}_k \rnb_1 H \left( \frac{\nu_{\mathbf{j},k} }{\lnb \nu_{\mathbf{j},k} \rnb_1 } \right) + \lnb \nu_{\mathbf{j},k} \rnb_1  \lnb \tilde{\eta}_k \rnb_1 H \left( \frac{\tilde{\eta}_k}{\lnb \tilde{\eta}_k \rnb_1 } \right)\\
&= \lnb \tilde{\eta}_k \rnb_1 H(\nu_{\mathbf{j},k}) + \lnb \nu_{\mathbf{j},k} \rnb_1 H(\tilde{\eta}_k).
\end{align*}
Therefore
\begin{align}
H(\mu * \tilde{\eta}_k) & =  H \left( \sum_{\mathbf{j} \in \left( \Z / k\Z \right)^d} \nu_{\mathbf{j},k} * \tilde{\eta}_k \right) \nonumber \\
& \geq  \sum_{\mathbf{j} \in \left( \Z / k\Z \right)^d} H \left(  \nu_{\mathbf{j},k} * \tilde{\eta}_k \right) \label{eq:entropy_tail_concave}\\
& \geq  \lnb \tilde{\eta}_k \rnb_1 H ( \mathbf{p}) + H(\tilde{\eta}_k) -  d \lnb \tilde{\eta}_k \rnb_1 \log k , \nonumber
\end{align}
where in \eqref{eq:entropy_tail_concave} we apply Lemma \ref{lemm:concave}. 

We wish to apply Lemma \ref{lemm:concave} again to sum over $k$. To do this we simply need to show that $\sum_{k=N}^{\infty} H(\mu * \tilde{\eta}_k)$ and $H(\sum_{k=N}^{\infty} \mu * \tilde{\eta}_k)$ both tend to zero as $N \to \infty$. In what follows, $c_1, c_2, \dots$ are positive constants, which depend only on $d$ and $R/r$. Note that we have
\begin{equation*}
\lnb \tilde{\eta}_k \rnb_1 \leq c_1 e^{- c_2 k^2}
\end{equation*}
and that the density function of $\tilde{\eta}_k$ is either $0$ or between $ \frac{c_3}{r} e^{- c_4 k^2}$ and $ \frac{c_5}{r} e^{- c_6 k^2}$. Also note that
\begin{equation*}
H(\tilde{\eta}_k) \leq H(\mu * \tilde{\eta}_k) \leq  \lnb \tilde{\eta}_k \rnb_1 H(\mu ) + H(\tilde{\eta}_k)
\end{equation*}
and so
\begin{equation*}
\left| H(\mu * \tilde{\eta}_k) \right| \leq c_7 e^{- c_8 k^2} \left( | \log r| + H(\mu)\right).
\end{equation*}
This means $\sum_{k=N}^{\infty} H(\mu * \tilde{\eta}_k) \to 0$. By our estimates on the density functions of $\tilde{\eta}_k$ we also have 
\begin{equation*}
\left| H(\sum_{k=N}^{\infty} \tilde{\eta}_k) \right| \leq c_9 e^{- c_{10} N^2}\left( | \log r| + 1 \right)
\end{equation*}
and so $H(\sum_{k=N}^{\infty} \mu * \tilde{\eta}_k) \to 0$. 

We then apply Lemma \ref{lemm:concave} to get
\begin{align}
H(\mu * \eta_{r^2}) & = H \left( \sum_{k=2}^{\infty} \mu*  \tilde{\eta}_k \right) \nonumber \\
& \geq \sum_{k=2}^{\infty} H \left(  \mu* \tilde{\eta}_k \right) \nonumber\\
&\geq H ( \mathbf{p}) + \sum_{k=2}^{\infty} H(\tilde{\eta}_k) - d \sum_{k=2}^{\infty} \lnb \tilde{\eta}_k \rnb_1 \log k. \label{eq:entropy_tail_temp_thing}
\end{align}

Recall that we have
\begin{equation*}
\lnb \tilde{\eta}_k \rnb_1 \leq c_{1} e^{- c_{2} k^2}
\end{equation*}
and so
\begin{equation}
H \left( \left( \lnb \tilde{\eta}_k \rnb_1 \right)_{k=2}^{\infty} \right) \leq c_{11} \label{eq:c_3_entropy_thing}
\end{equation}
and
\begin{equation*}
d \sum_{k=2}^{\infty} \lnb \tilde{\eta}_k \rnb_1 \log k  \leq c_{12}.
\end{equation*}

Applying Lemma \ref{lemm:max_inc} and \eqref{eq:c_3_entropy_thing}, we have
\begin{align*}
\MoveEqLeft d \log r + c_{13} = H(\eta_{r^2}) = H \left( \sum_{k=2}^{\infty} \tilde{\eta}_k \right)\\
&\leq \sum_{k=2}^{\infty} H \left(  \tilde{\eta}_k \right) + H \left( \left( \lnb \tilde{\eta}_k \rnb_1 \right)_{k=2}^{\infty} \right) \leq \sum_{k=2}^{\infty} H \left(  \tilde{\eta}_k \right) + c_{14}.
\end{align*}
Substituting this estimate for $\sum_{k=2}^{\infty} H \left(  \tilde{\eta}_k \right)$ into \eqref{eq:entropy_tail_temp_thing} gives the required result.
\end{proof}

\subsection{Proof of Lemma \ref{lemm:sum_expression_lemma}}

In order for Lemma \ref{lemm:initial_gap} to be useful it is necessary to show that if $I_1, I_2, \dots, I_n$ are disjoint intervals contained in $(0, 1]$ then there exists some $\nu$ such that $\mu_F = \nu * \mu_F^{I_1} * \mu_F^{I_2} * \dots \mu_F^{I_n}$. To do this it suffices to prove Lemma \ref{lemm:sum_expression_lemma}. Indeed we can then take $\nu = \mu_F^{(0,1]\backslash (I_1 \cup I_2 \cup \dots \cup I_n)}$.

\begin{proof}[Proof of Lemma \ref{lemm:sum_expression_lemma}]
For $k$ in $\N$ let $Y_k$ be defined by
\begin{equation*}
Y_k := \sum_{i=0}^{k-1} \lambda^i U^i X_i
\end{equation*}
and let $\mu_k$ be the law of $Y_k$. It is clear that $\mu_k$ satisfies
\begin{equation}
\mu_{k+1} = \sum_{i=1}^{n} p_i \mu_k \circ S_i^{-1}. \label{eq:step_eq}
\end{equation}
Let $\mu$ be the law of $Y$. Clearly we have that $Y_k \to Y$ almost surely and so $\mu_{k}$ tends to $\mu$ weakly. Taking the weak limit of both sides of \eqref{eq:step_eq} gives
\begin{equation*}
\mu = \sum_{i=1}^{n} p_i \mu \circ S_i^{-1}. 
\end{equation*}
Therefore by the uniqueness of $\mu_F$ we get that $\mu=\mu_F$ as required.
\end{proof}

\section{Proof of the main theorem} \label{section:main_proof}

We follow the strategy outlined in Section \ref{section:outing_of_proof}. To implement this we make the following definition.

\begin{defi} \label{defi:admissible_interval}
Given some $r \in \left(0, \frac{1}{10} \right)$ and iterated function system $F$ on $\R^d$ we say that an interval $I \subset (0, \infty)$ is \emph{$\alpha$-admissible at scale $r$} if for all $t$ with
\begin{equation*}
t \in \left[ \exp \left( -\left( \log \log r^{-1} \right) ^{10} \right)r,  \exp \left( \left( \log \log r^{-1} \right) ^{10} \right)r\right],
\end{equation*}
we have
\begin{equation*}
\left. \frac{\partial}{\partial y} H(\mu_F^{I}*\eta_y) \right|_{y=t^2} \leq \alpha t^{-2} .
\end{equation*}
\end{defi}

Recall that $\mu_F^I$ is as defined in Definition \ref{defi:pieces}. This definition is designed in such a way that if $I_1$ and $I_2$ is a pair of disjoint admissible intervals, then we can apply Theorem \ref{prop:entropy_to_detail} for the measure $\mu_F^{I_1 \cup I_2} = \mu_F^{I_1} * \mu_F^{I_2}$ to obtain estimates for $s_t(\mu_F^{I_1 \cup I_2})$ at a range of scales $t$ in a suitable range around $r$. Moreover, these estimates are suitable so that we can apply Theorem \ref{theo:many_conv} for $\mu_F^{I_1 \cup I_2}$ in the role of one of the measures. If we have many admissible intervals we get an improved estimate for $s_r(\mu_F)$ via Theorem \ref{theo:many_conv}.

We formalize the result of these ideas in the following statement. The detail of its proof is given in Section \ref{section:detail_of_admissible}.

\begin{prop} \label{prop:many_admissible_to_detail}
Let $\alpha, K >0$ and let $d \in \N$. Suppose that $\alpha < \frac{1}{8}\left(1 + \frac{1}{2K^2} \right)^{-1}$. Then there exists some constant $c$ such that the following is true.
	
Let $F$ be an iterated function system on $\R^d$ with uniform contraction ratio and uniform rotation. Suppose that $r <c$ and  $n\in \N$ is even with
\begin{equation}
n \leq 10 \frac{\log \log r^{-1}}{\log \left( \frac{1}{8} \left(1 + \frac{1}{2K^2} \right)^{-1} \alpha^{-1} \right)} \label{eq:max_num_intervals}
\end{equation}
and that $I_1, I_2, \dots, I_n$ are disjoint $\alpha$-admissible intervals at scale $r$ contained in $(0,1)$. Then we have
\begin{equation}
s_r(\mu_F) \leq  \frac{1}{4} Q(d) \left( 8\left(1 + \frac{1}{2K^2} \right) \alpha \right)^{\frac{n}{2}} . \label{admissibledetail}
\end{equation}
\end{prop}

Our next goal is to find suitably many disjoint admissible intervals at a given scale $r$. This is done using Lemma \ref{lemm:initial_gap} in Section \ref{section:find} where we prove the following statement.

\begin{lemm} \label{lemm:find}
Suppose that $F$ is an iterated function system with uniform rotation and uniform contraction ratio $\lambda$. Let $M>M_F$, $\alpha \in (0, \frac{1}{8})$ and suppose that $P>1$ and satisfies
\begin{equation}
d \log M - h_F < 2 \alpha \lambda^2 ( \log M - P \log \lambda^{-1}). \label{eq:inequality_in_find}
\end{equation}
Then there exists some $c >0$ such that for every $r>0$ sufficiently small there are at least
\begin{equation*}
\frac{1}{\log \frac{\log M}{(P-1) \log \lambda^{-1}}} \log \log r^{-1} - c \log \log \log r^{-1}
\end{equation*}
disjoint $\alpha$-admissible intervals at scale $r$ all of which are contained in $(0,1]$ . 
\end{lemm}

It is worth pointing out that we always have $h_F \leq d \log M_F$ and $h_F$ can be arbitrarily close to this upper limit. This means that \eqref{eq:inequality_in_find} can be satisfied for any given value of $\alpha$ and $P$ provided $h_F$ is sufficiently close to $d \log M_F$ and $M$ is sufficiently close to $M_F$.

In order to apply Lemma \ref{lemm:suff_abs_cont}, we wish to show that $s_r(\mu_F) \leq \left( \log r^{-1} \right)^{-\beta}$ for some $\beta > 1$ for all sufficiently small $r$. Since we may take $K$ arbitrarily large in Proposition \ref{prop:many_admissible_to_detail}, it suffices to show that there is some $\beta > 1$ such that for all sufficiently small $r>0$, we can find at least $\beta \frac{2 \log \log r^{-1}}{\log 1/(8\alpha)}$ disjoint admissible intervals. In Section \ref{section:subsection_main_proof} we use Lemma \ref{lemm:find} and a careful choice of $\alpha$ and $P$ to do this.

The condition \eqref{eq:max_num_intervals} is unimportant because if we have more than this many admissible intervals, then it turns out that taking $n$ to be the greatest even number less than
\begin{equation*}
10 \frac{\log \log r^{-1}}{\log \left( \frac{1}{8} \left(1 + \frac{1}{2K^2} \right)^{-1} \alpha^{-1} \right)} 
\end{equation*}
gives a sufficiently strong bound on detail to prove absolute continuity.

\subsection{Detail of the convolution of many admissible pieces} \label{section:detail_of_admissible}

In this subsection, we prove Proposition \ref{prop:many_admissible_to_detail}.

\begin{proof}[Proof of Proposition \ref{prop:many_admissible_to_detail}]
Throughout this proof, let $c_1, c_2, ...$ denote constants depending only on $\alpha, K$ and $d$. The idea is to use Theorem \ref{theo:many_conv} and Proposition \ref{prop:entropy_to_detail}.

First note that by applying Proposition \ref{prop:entropy_to_detail} with $u=v=\frac{t^2}{2}$ we know that for all 
\begin{equation*} 
t \in \left[ \sqrt{2} \exp \left( -\left( \log \log r^{-1} \right) ^{10} \right)r,  \sqrt{2} \exp \left( \left( \log \log r^{-1} \right) ^{10} \right)r\right]
\end{equation*}
and for $i = 1, 2, ... , \frac{n}{2}$ we have
\begin{align*}
s_t(\mu_F^{I_{2i-1}} * \mu_F^{I_{2i}}) & \leq t^2 Q(d) \sqrt{ \left. \frac{\partial}{\partial y} H \left( \mu_F^{I_{2i-1}} * \eta_y \right) \right|_{y = \frac{t^2}{2}}  \left. \frac{\partial}{\partial y} H \left( \mu_F^{I_{2i}} * \eta_y \right) \right|_{y = \frac{t^2}{2}} } \\
& \leq 2 Q(d) \alpha.
\end{align*}
We now wish to apply Theorem \ref{theo:many_conv} for the measures $\mu_F^{I_{2i-1}\cup I_{2i}}$ for $i=1, 2, \dots \frac{n}{2}$ with $\alpha_1 = \alpha_2 = \dots = \alpha_{n/2} = 2 Q(d) \alpha $. To do this we simply need to check that
\begin{equation*}
\left[ 2^{-\frac{m}{2}} r, K^m \alpha_1^{-m 2^m} r \right] \subset \left[ \sqrt{2} \exp \left( -\left( \log \log r^{-1} \right) ^{10} \right)r,  \sqrt{2}\exp \left( \left( \log \log r^{-1} \right) ^{10} \right)r\right]
\end{equation*}
where $m = \frac{\log (n/ 2)}{\log (3/2)}$. We note that
\begin{equation*}
m \leq \frac{1}{\log (3/2)} \log \log \log r^{-1} + c_1
\end{equation*}
and so for all sufficiently small $r$, we have
\begin{equation*}
2^{-\frac{m}{2}} r \geq \sqrt{2} \exp \left( -\left( \log \log r^{-1} \right) ^{10} \right)r .
\end{equation*}
For the other side, note that
\begin{equation*}
K^m \alpha_1^{-m 2^m}  \leq \exp \left( c_2 \left( \log \log \log r^{-1} \right) \left( \log \log r^{-1} \right)^{\frac{\log 2}{\log (3/2)}} +c_3 \right) .
\end{equation*}
Noting that $\frac{\log 2}{\log (3/2)}<10$, for all sufficiently small $r$ we have
\begin{equation*}
K^m \alpha_1^{-m 2^m} r \leq \exp \left( \left( \log \log r^{-1} \right) ^{10} \right)r .
\end{equation*}
Therefore, the conditions of Theorem \ref{theo:many_conv} are satisfied and so
\begin{align*}
\MoveEqLeft s_r(\mu_F^{I_1} * \mu_F^{I_2} * \dots * \mu_F^{I_n}) \\
& \leq \left( 2 Q(d) \alpha \right) ^{\frac{n}{2}} \left( \frac{4}{Q(d)} \left( 1 + \frac{1}{2K^2} \right) \right) ^{\frac{n}{2}-1} \\
& \leq \frac{1}{4} Q(d) \left( 8\left(1 + \frac{1}{2K^2} \right) \alpha \right)^{\frac{n}{2}} .
\end{align*}
We conclude the proof by noting that by Proposition \ref{prop:no_detail_increase}
\begin{equation*}
s_r(\mu_F) \leq s_r(\mu_F^{I_1} * \mu_F^{I_2} * \dots * \mu_F^{I_n}) . \qedhere
\end{equation*}
\end{proof}

\subsection{Finding admissible intervals} \label{section:find}

In this subsection, we prove Lemma \ref{lemm:find}. The main ingredient in the proof of Lemma \ref{lemm:find} is the following lemma.

\begin{lemm} \label{lemm:sub_find}
Let $F$ be an iterated function system with uniform rotation and uniform contraction ratio $\lambda$. Let $\alpha, r > 0$, $n \in \Z_{\geq 0}$ and $k \in \Z$. Suppose that
\begin{equation*}
\frac{\partial}{\partial y} H(\mu_F^{(\lambda^n, 1]} * \eta_y) \leq \frac{1}{y} \lambda^2 \alpha
\end{equation*}
for some $y \in \left( \lambda^{2k+2}, \lambda^{2k} \right]$. Then the interval
\begin{equation}
I = \left(r\lambda^{ n - k + b(r)}, r\lambda^{-k   - b(r) }\right] \label{eq:def_i_sub}
\end{equation}
is $\alpha$-admissible at scale $r$. Here $b=b(r)$ is an error term defined by
\begin{equation*}
b :=  \frac{1}{\log \lambda^{-1}} \left( \log \log r^{-1} \right) ^{10} +10.
\end{equation*}
\end{lemm}

We first prove Lemma \ref{lemm:sub_find} and then proceed with the proof of Lemma \ref{lemm:find}. To prove this, we need a few more facts about entropy. It is well known that for any absolutely continuous random variable $X$ taking values in $\R^d$ and any bijective linear map $A:\R^d\to\R^d$ we have
\begin{equation*}
H(AX) = H(X) + \log |\det A|.
\end{equation*}
It follows that
\begin{equation*}
H(\mu_F^{(\lambda^k, \lambda^{\ell}]} * \eta_{t^2}) = H(\mu_F^{(\lambda^{k-\ell}, 1]} * \eta_{\lambda^{-2\ell} t^2}) + d \ell \log \lambda
\end{equation*}
and also
\begin{equation}
\left. \frac{\partial}{\partial y} H(\mu_F^{(\lambda^k, \lambda^{\ell}]} * \eta_y) \right|_{y = t^2} = \lambda^{-2\ell} \left. \frac{\partial}{\partial y}  H(\mu_F^{(\lambda^{k-\ell}, 1]}*\eta_y) \right|_{y = \lambda^{-2\ell} t^2}. \label{eq:rescale1}
\end{equation}

We also have the following.
\begin{prop} \label{prop:madiman}
Let $X_1$, $X_2$ and $X_3$ be independent absolutely continuous random variables with finite entropy. Then,
\begin{equation*}
H(X_1 + X_2 + X_3) + H(X_1) \leq H(X_1 + X_2) + H(X_1 + X_3) .
\end{equation*}
\end{prop}
\begin{proof}
This is proven in \cite[Theorem 3.1]{Madiman_Kontoyiannis_2014}.
\end{proof}

\begin{coro} \label{coro:decrease_ny}
Let $\mu$ and $\nu$ be measures on $\R^d$ with finite variance and let $y > 0$. Then
\begin{equation}
\frac{\partial}{\partial y} H(\mu * \nu * \eta_y) \leq \frac{\partial}{\partial y} H(\mu * \eta_y).
\end{equation}

\end{coro}
\begin{proof}
Let $\varepsilon > 0$. Then using Proposition \ref{prop:madiman} with $X_1, X_2$ and $X_3$ having laws $\mu* \eta_{y}$, $\eta_{\varepsilon}$ and $\nu$ respectively we get
\begin{equation*}
H(\mu * \nu * \eta_{y} * \eta_{\varepsilon}) - H(\mu * \nu * \eta_{y}) \leq H(\mu * \eta_{y} * \eta_{\varepsilon}) - H(\mu * \eta_{y}).
\end{equation*}
The result follows by taking the limit $\varepsilon \to 0$.
\end{proof}

An immediate consequence of Corollary \ref{coro:decrease_ny} is that the function $y \mapsto \frac{\partial}{\partial y} H(\mu * \eta_y)$ is non-increasing and if $I_1 \subset I_2$ then
\begin{equation}
\frac{\partial}{\partial y} H(\mu_F^{I_2} * \eta_y) \leq \frac{\partial}{\partial y} H(\mu_F^{I_1} * \eta_y) \label{eq:entropy_smaller_interval_contained}.
\end{equation}
In particular this means that if $I_1$ is $\alpha$-admissible at scale $r$ for some $\alpha$ and $r$ then so is $I_2$. This is important both for proving Lemma \ref{lemm:sub_find} and for showing that Lemma \ref{lemm:find} follows from Lemma \ref{lemm:sub_find}. We are now ready to prove Lemma \ref{lemm:sub_find}.

\begin{proof}[Proof of Lemma \ref{lemm:sub_find}]
To prove this, suppose that
\begin{equation*}
t \in \left[ \exp \left( -\left( \log \log r^{-1} \right) ^{10} \right)r,  \exp \left( \left( \log \log r^{-1} \right) ^{10} \right)r\right].
\end{equation*}
We wish to show that $\left. \frac{\partial}{\partial y} H(\mu_F^I* \eta_y) \right|_{y=t^2} \leq \alpha t^{-2}$, where $I$ is defined in \eqref{eq:def_i_sub}. Choose $\tilde{t} \in (\lambda^{k+1}, \lambda^k]$ such that 
\begin{equation}
\left. \frac{\partial}{\partial y} H(\mu_F^{(\lambda^n, 1]} * \eta_y) \right|_{y = \tilde{t}^2} \leq \alpha \lambda^2 \tilde{t}^{-2}\label{eq:tilde_t_def}
\end{equation}
 and choose $\tilde{k} \in \mathbb{Z}$ such that 
\begin{equation*}
\lambda^{\tilde{k} + 1} \tilde{t} \leq t \leq \lambda ^ {\tilde{k}} \tilde{t}.
\end{equation*}
We then have
\begin{align}
\left. \frac{\partial}{\partial y} H \left(\mu_F^{(\lambda^{n + \tilde{k} +1},\lambda^{\tilde{k} +1}]} * \eta_y \right) \right|_{y=t^2} & \leq \left. \frac{\partial}{\partial y} H \left(\mu_F^{(\lambda^{n + \tilde{k} +1},\lambda^{\tilde{k} +1}]} * \eta_y \right) \right|_{y=\lambda^{2 \tilde{k} + 2} \tilde{t}^2} \label{eq:n_ydec}\\
& = \lambda^{-2 \tilde{k} -2} \left. \frac{\partial}{\partial y} H \left(\mu_F^{(\lambda^n,1]} * \eta_y \right) \right|_{y=\tilde{t}^2} \label{eq:rescale2} \\
& \leq \lambda^{-2 \tilde{k}-2} \cdot \lambda^2 \alpha \tilde{t} ^{-2} \label{eq:uses_ttilde_in_stilde} \\
& \leq \alpha t ^{-2}. \nonumber
\end{align}
Where \eqref{eq:n_ydec} follows from Corollary \ref{coro:decrease_ny}, \eqref{eq:rescale2} follows from \eqref{eq:rescale1} and \eqref{eq:uses_ttilde_in_stilde} follows from \eqref{eq:tilde_t_def}.

Note that $\left(\lambda^{n + \tilde{k} +1}, \lambda^{\tilde{k}+1} \right] \subset I$ hence by \eqref{eq:entropy_smaller_interval_contained} we have
\begin{equation*}
\left.\frac{\partial}{\partial y} H \left( \mu_F^{I} * \eta_y \right)\right|_{y=t^2} \leq \alpha t^2
\end{equation*}
as required.
\end{proof}

We can use Lemma \ref{lemm:sub_find} and Lemma \ref{lemm:initial_gap} to show that some specific intervals are $\alpha$-admissible at scale $r$. We prove the following.

\begin{lemm} \label{lemm:interval_is_admissible_b}
Suppose that $F$ is an iterated function system with uniform contraction ratio $\lambda$ and uniform rotation and that $M>M_F$. Let $\alpha \in (0,1)$. Suppose further that there is some constant $P>1$ such that
\begin{equation}
d \log M - h_F < 2 \alpha \lambda^2 (\log M - P \log \lambda^{-1}). \label{eq:new_condition_to_find_interval}
\end{equation}
Then for all sufficiently large $n \in \N$ and all $r \in (0,\frac{1}{4})$ the interval

\begin{equation*}
I = \left( r\lambda^{k_1},r \lambda^{k_2} \right]
\end{equation*}
is $\alpha$-admissible at scale $r$.

Here $b = b(r)$ be defined by
\begin{equation*}
b :=  \frac{1}{\log \lambda^{-1}} \left( \log \log r^{-1} \right) ^{10} +10 ,
\end{equation*}
$k_1$ is defined by
\begin{equation*}
k_1 :=  -(P-1) n  +  b(r)
\end{equation*}
and $k_2$ is defined by
\begin{equation*}
k_2 :=  - n \frac{\log M}{\log \lambda^{-1}}  - b(r).
\end{equation*}

\end{lemm}

\begin{proof}
Suppose for contradiction that this is not true. Recall that if $I_1 \subset I_2$ and $I_1$ is $\alpha$-admissible at scale $r$ then $I_2$ is $\alpha$-admissible at scale $r$. Therefore by Lemma \ref{lemm:sub_find} we have that there cannot exist $k \in \Z_{\geq 0}$ and $y \in \left( \lambda^{2k+2}, \lambda^{2k}\right]$ such that
\begin{equation*}
\frac{\partial}{\partial y} H(\mu_F^{(\lambda^n, 1]} * \eta_y) \leq \frac{1}{y} \lambda^2 \alpha.
\end{equation*}
and
\begin{equation}
k_1 =  -(P-1)n  + b(r) \geq  + n - k + b(r) \label{eq:temp_admiss_lemma_1}
\end{equation}
and
\begin{equation}
k_2 =  - \frac{\log M}{\log \lambda^{-1}} n - b(r) \geq  - k - b(r). \label{eq:temp_admiss_lemma_2}
\end{equation}

Note that \eqref{eq:temp_admiss_lemma_1} is equivalent to $k \geq Pn - c$ and \eqref{eq:temp_admiss_lemma_2} is equivalent to $k \leq \frac{\log M}{\log \lambda^{-1}}n$. In particular, noting that $\lambda^{2 \frac{\log M}{\log \lambda^{-1}} n} = M^{-2n}$, this means that 
we have
\begin{equation*}
\frac{\partial}{\partial y} H(\mu_F^{(\lambda^n, 1]} * \eta_y) > \frac{1}{y} \lambda^2 \alpha
\end{equation*}
for all $y$ such that
\begin{equation*}
y \in \left( M^{-2n}, \lambda^{2 Pn } \right].
\end{equation*}

In particular, provided $n$ is sufficiently large, by integrating we get
\begin{align}
\MoveEqLeft H\left(\mu_F^{(\lambda^n, 1]} * \eta_1 \right) - H\left(\mu_F^{(\lambda^n, 1]} * \eta_{M^{-2n}} \right) \nonumber\\
& \geq H\left(\mu_F^{(\lambda^n, 1]} * \eta_{\lambda^{2 Pn}} \right) - H\left(\mu_F^{(\lambda^n, 1]} * \eta_{M^{-2n}} \right) \label{eq:ent_monotone_subset}\\
& = \int_{M^{-2n}}^{\lambda^{2 Pn}} \frac{\partial}{\partial y} H\left(\mu_F^{(\lambda^n, 1]} * \eta_y \right) \, dy \nonumber \\
& \geq \int_{M^{-2n}}^{\lambda^{2 Pn}} \frac{1}{y} \alpha \lambda^2 \, dy  \nonumber \\
& = 2 n \alpha \lambda^2 \left(\log M - P \log \lambda^{-1} \right)  \nonumber
\end{align}
with \eqref{eq:ent_monotone_subset} following from Lemma \ref{lemm:ent_increase}. This contradicts Lemma \ref{lemm:initial_gap}.
\end{proof}

We are now ready to prove Lemma \ref{lemm:find}.

\begin{proof}[Proof of Lemma \ref{lemm:find}]
Throughout this proof $E_1, E_2, \dots$ denote error terms which may be bounded by $0 \leq E_i \leq c_i \left( \log \log r^{-1} \right)^{c_i}$ for some positive constants $c_1, c_2, \dots$ which depend only on $\alpha$, $F$, $P$ and $M$. 
Let $c'$ take the role of $c$ in Lemma \ref{lemm:interval_is_admissible_b} and choose $N$ large enough that Lemma \ref{lemm:interval_is_admissible_b} holds for all $n \geq N$.

We wish to choose some $j_{\text{max}}$ and some $N= n_0 < n_1 < n_2 < \dots < n_{j_{\text{max}}}$ such that if we let
\begin{equation*}
k_1^{(j)} = \frac{\log r^{-1}}{\log \lambda^{-1}} - (P-1) n_j + c' +b
\end{equation*}
and
\begin{equation*}
k_2^{(j)} = \frac{\log r^{-1}}{\log \lambda^{-1}} - \frac{\log M}{\log \lambda^{-1}} n_j - b
\end{equation*}
and
\begin{equation*}
I_j = \left( \lambda^{k_1^{(j)}}, \lambda^{k_2^{(j)}} \right]
\end{equation*}
then each of $I_0, I_1, \dots, I_{j_{\text{max}}}$ are disjoint subsets of $(0,1]$. Note that by Lemma \ref{lemm:interval_is_admissible_b}, each of the $I_j$ are $\alpha$-admissible at scale $r$. In order for the intervals to be disjoint it is sufficient to have $k_2^{(j)} \geq k_1^{(j+1)}$ for $j = 0, 1, \dots, j_{\text{max}}-1$. This is equivalent to
\begin{equation*}
\frac{\log r^{-1}}{\log \lambda^{-1}} - \frac{\log M}{\log \lambda^{-1}} n_j - b \geq \frac{\log r^{-1}}{\log \lambda^{-1}} - (P-1) n_{j+1} + c' +b
\end{equation*}
which becomes
\begin{equation}
n_{j+1} \geq \frac{\log M}{(P-1) \log \lambda^{-1}} n_j + E_1. \label{eq:in_find_proof}
\end{equation}

Note that by the hypothesis of Lemma \ref{lemm:find} we have $\log M \geq P \log \lambda^{-1} > (P-1) \log \lambda^{-1}$ and so $\frac{\log M}{(P-1) \log \lambda^{-1}} > 1$.

We achieve \eqref{eq:in_find_proof} by taking $n_{j+1} = \ceil{\frac{\log M}{(P-1) \log \lambda^{-1}} n_j + E_1}$. Note that this gives $n_{j+1} \leq \frac{\log M}{(P-1) \log \lambda^{-1}} n_j + E_2$ which can be rewritten as
\begin{equation*}
n_{j+1} + \frac{1}{\frac{\log M}{(P-1) \log \lambda^{-1}}-1} E_2 \leq \frac{\log M}{(P-1) \log \lambda^{-1}} \left( n_j + \frac{1}{\frac{\log M}{(P-1) \log \lambda^{-1}}-1} E_2\right)
\end{equation*}
which gives
\begin{equation}
n_j \leq \left( \frac{\log M}{(P-1) \log \lambda^{-1}} \right)^j(n_0+E_3). \label{eq:n_j_upper_bound}
\end{equation}
Noting that $n_0 = N = E_4$ we get
\begin{equation*}
n_j \leq \left( \frac{\log M}{(P-1) \log \lambda^{-1}} \right)^j E_5.
\end{equation*}

We also need to ensure that all of the intervals $I_0, I_1, \dots, I_{j_{\text{max}}}$ are contained in $(0,1]$. For this it is sufficient to show that
\begin{equation*}
\frac{\log r^{-1}}{\log \lambda^{-1}} - \frac{\log M}{\log \lambda^{-1}} n_{j_{\text{max}}} - E_6 \geq 0.
\end{equation*}
By \eqref{eq:n_j_upper_bound} it is sufficient to have
\begin{equation*}
\left( \frac{\log M}{(P-1) \log \lambda^{-1}} \right)^{j_{\text{max}}} E_7 \leq \log r^{-1}
\end{equation*}
which can be achieved with
\begin{equation*}
j_{\text{max}} = \ceil{\frac{1}{\log \frac{\log M}{(P-1) \log \lambda^{-1}}} \log \log r^{-1} - c \log \log \log r^{-1}}
\end{equation*}
for some constant $c$ depending only on $\alpha$, $F$ and $M$ for all sufficiently small $r$ as required. In particular this gives
\begin{equation*}
j_{\text{max}} \geq \frac{1}{\log \frac{\log M}{(P-1) \log \lambda^{-1}}} \log \log r^{-1} - c \log \log \log r^{-1}
\end{equation*}
as required.
\end{proof}

\subsection{Proof of the main theorem} \label{section:subsection_main_proof}

We are now ready to prove Theorem \ref{theo:main_result}.

\begin{proof}[Proof of Theorem \ref{theo:main_result}]
The idea is to use Proposition \ref{prop:many_admissible_to_detail} and Lemma \ref{lemm:find} to show that the detail around a scale decreases fast enough for us to be able to apply Lemma \ref{lemm:suff_abs_cont}.

Let $M > M_F$ and throughout this proof let $c_1, c_2, \dots$ denote constants that depend only on $M$, $F$, $P$ and $\alpha$. Note that by Lemma \ref{lemm:find} given any $M > M_F$ for all sufficiently small $r$ there are at least
\begin{equation*}
\frac{1}{\log A} \log \log r^{-1} - c_1 \log \log \log r^{-1}
\end{equation*}
disjoint admissible intervals contained in $(0,1]$ where
\begin{equation*}
A = \frac{\log M}{(P-1) \log \lambda^{-1}} 
\end{equation*}

By Proposition  \ref{prop:many_admissible_to_detail}, we have that
\begin{equation*}
s_r(\mu_F) \leq c_2 \left( 8 \left( 1 + \frac{1}{2K^2} \right) \alpha \right)^{n/2},
\end{equation*}
where  $n$ is the largest even number which is less than both $\frac{1}{\log A} \log \,\log \, r^{-1} - c_1 \log \, \log \, \log \, r^{-1}$ and $10 \frac{\log \log r^{-1}}{\log \left( \frac{1}{8} \left(1 + \frac{1}{2K^2} \right)^{-1} \alpha^{-1} \right)}$.

If
\begin{equation*}
\frac{1}{\log A} \geq \frac{10}{\log \left( \frac{1}{8} \left(1 + \frac{1}{2K^2} \right)^{-1} \alpha^{-1} \right)}
\end{equation*}
then 
\begin{equation*}
n \geq \frac{10}{\log \left( \frac{1}{8} \left(1 + \frac{1}{2K^2} \right)^{-1} \alpha^{-1} \right)} - c_3 \log \, \log \, \log \, r^{-1}
\end{equation*}
and so
\begin{align*}
s_r(\mu_F) & \leq c_2 \exp \left( - 5 \log \log r^{-1} + c_4 \log \log \log r^{-1} \right)\\
& = c_2 \left( \log r^{-1} \right) ^{-5} \left(\log \log r^{-1} \right) ^{c_4}.
\end{align*}
By Lemma \ref{lemm:suff_abs_cont} it follows that $\mu_F$ is absolutely continuous. 

If instead
\begin{equation*}
\frac{1}{\log A} < \frac{10}{\log \left( \frac{1}{8} \left(1 + \frac{1}{2K^2} \right)^{-1} \alpha^{-1} \right)}
\end{equation*}
then we get 
\begin{equation*}
n \geq \frac{1}{\log A} \log \log r^{-1} - c_3 \log \, \log \, \log \, r^{-1}.
\end{equation*}
This gives
\begin{align*}
s_r(\mu_F) & \leq c_2 \left( 8 \left( 1 + \frac{1}{2K^2} \right) \alpha \right) ^ {\frac{1}{2 \log A} \log \log r^{-1} - c_5 \log \, \log \, \log \, r^{-1}}\\
& = c_2\exp \left( - \frac{ \log \left( 8 \left( 1 + \frac{1}{2 K^2} \right) \alpha \right) ^{-1}}{2 \log A} \log \log r^{-1} + c_6 \log \log \log r^{-1} \right)\\
& = c_2 \left( \log r^{-1} \right)^{-\beta} \left( \log \log r^{-1} \right)^{c_7}
\end{align*}
where $\beta = \frac{ \log \left( 8 \left( 1 + \frac{1}{2 K^2} \right) \alpha \right) ^{-1}}{2 \log A}$.

By Lemma \ref{lemm:suff_abs_cont} for $\mu_F$ to be absolutely continuous it is sufficient to have $\beta > 1$. For this it is sufficient to show that
\begin{equation*}
\left( 8 \left( 1 + \frac{1}{2 K^2} \right) \alpha \right) ^{-1} > A^2.
\end{equation*}

Since we can choose $K$ to be arbitrarily large and $M$ to be arbitrarily close to $M_F$ it is sufficient to have
\begin{equation}
\frac{1}{8 \alpha} > \tilde{A}^2 \label{eq:sufficient_inequality}
\end{equation}
where
\begin{equation*}
\tilde{A} = \frac{\log M_F}{(P-1) \log \lambda^{-1}}.
\end{equation*}

Also by choosing $M$ sufficiently close to $M_F$ our condition on $P$ becomes
\begin{equation*}
d \log M_F - h_F < 2 \alpha \lambda^2 (\log M_F - P \log \lambda^{-1})
\end{equation*}
which may be written as
\begin{equation*}
P < \frac{d \log M_F - h_F - 2 \alpha \lambda^2 \log M_F}{2 \alpha \lambda^2 \log \lambda^{-1}}.
\end{equation*}
By choosing $P$ arbitrarily close to this upper bound and taking the square root of both sides \eqref{eq:sufficient_inequality} becomes
\begin{equation*}
\frac{1}{\sqrt{8\alpha}} > \frac{2 \alpha \lambda^2 \log M_F}{h_F - 2 \alpha \lambda^2 \log \lambda^{-1} - (d - 2 \alpha \lambda^2) \log M_F}
\end{equation*}
which can be rewritten as
\begin{equation}
h_F - 2 \alpha \lambda^2 \log \lambda^{-1} - (d - 2 \alpha \lambda^2) \log M_F > \sqrt{8 \alpha} (2 \alpha \lambda^2 \log M_F). \label{eq:condition_for_alpha}
\end{equation}

We now substitute in $\alpha = \frac{1}{18} \left( \frac{ \log M_F - \log \lambda^{-1}}{\log M_F} \right)^2$ (it is easy to check by differentiating \eqref{eq:condition_for_alpha} that this is the optimal choice for $\alpha$). The inequality becomes
\begin{align*}
\MoveEqLeft h_F - \frac{1}{9} \left( \frac{ \log M_F - \log \lambda^{-1}}{\log M_F} \right)^2 \lambda^2 \log \lambda^{-1} - \left(d - \frac{1}{9} \left( \frac{ \log M_F - \log \lambda^{-1}}{\log M_F} \right)^2 \lambda^2 \right) \log M_F  \\ & > \frac{2}{3} \left( \frac{ \log M_F - \log \lambda^{-1}}{\log M_F} \right) \left( \frac{1}{9} \left( \frac{ \log M_F - \log \lambda^{-1}}{\log M_F} \right)^2 \lambda^2 \log M_F \right).
\end{align*}
Multiplying both sides by $\left(\log M_F \right)^2$ gives
\begin{align*}
\MoveEqLeft h_F \left(\log M_F \right)^2 - \frac{1}{9} \left( \log M_F - \log \lambda^{-1} \right)^2 \lambda^2 \log \lambda^{-1}\\
\MoveEqLeft - \left(d \left(\log M_F \right)^2 - \frac{1}{9} \left( \log M_F - \log \lambda^{-1} \right)^2 \lambda^2 \right) \log M_F \\ & > \frac{2}{27} \left(  \log M_F - \log \lambda^{-1}\right) \left( \left(  \log M_F - \log \lambda^{-1} \right)^2 \lambda^2  \right).
\end{align*}
Rearranging reduces the inequality to
\begin{equation*}
(d \log M_F - h_F)(\log M_F)^2 < \frac{1}{27} \left( \log M_F - \log \lambda^{-1} \right) ^3 \lambda^2 
\end{equation*}
as required.

We now simply need to check that we have $P>1$. Since we choose $P$ arbitrarily close to $\frac{d \log M_F - h_F - 2 \alpha \lambda^2 \log M_F}{2 \alpha \lambda^2 \log \lambda^{-1}}$ it suffices to show that
\begin{equation*}
\frac{d \log M_F - h_F - 2 \alpha \lambda^2 \log M_F}{2 \alpha \lambda^2 \log \lambda^{-1}} > 1.
\end{equation*}
With our choice of $\alpha$ this becomes
\begin{equation*}
d \log M - h_F < \frac{1}{9} \left( \frac{\log M_F - \log \lambda^{-1}}{\log M_F} \right)^2 \lambda^2 (\log M - \log \lambda^{-1})
\end{equation*}
which may be rewritten as
\begin{equation*}
\left( d \log M - h_F \right) \left( \log M_F \right)^2 < \frac{1}{9} \left( \log M_F - \log \lambda^{-1} \right)^2 \lambda^2 (\log M - \log \lambda^{-1}).
\end{equation*}
Clearly this is satisfied under the conditions of Theorem \ref{theo:main_result} provided $M$ is sufficiently close to $M_F$ as required.
\end{proof}

\subsection{Proof of the result for Bernoulli convolutions} \label{section:main_proof_bernoulli}

We also wish to explain how Theorem \ref{theo:main_bernoulli} follows from Theorem \ref{theo:main_result}. First of all we use the following lemma to bound $M_F$.

\begin{lemm} \label{lemm:garsia} Let $\lambda$ be an algebraic number and denote by $d$ the number of its algebraic conjugates with modulus $1$. Then there is some constant $c_{\lambda}$ depending only on $\lambda$ such that whenever $p$ is a polynomial with degree $n$ and coefficients $-1, 0$ and $1$ such that $p(\lambda) \neq 0$ we have
\begin{equation*}
|p(\lambda)| > c_{\lambda} n^{-d} M_{\lambda} ^ {-n} .
\end{equation*}
\end{lemm}

\begin{proof}
This is proven in \cite[Lemma 1.51]{GARSIA_1962}.
\end{proof}

\begin{coro} \label{coro:garsia}
Let $F$ be an iterated function system such that  $\mu_F$ is a Bernoulli convolution with parameter $\lambda$. Then
\begin{equation*}
M_F \leq M_{\lambda} .
\end{equation*}
\end{coro}

\begin{proof}
If $x$ and $y$ are both in the support of $\sum_{i=0}^{n-1} \pm \lambda^i$ then clearly $x-y = 2 p(\lambda)$ for some polynomial $p$ of degree at most $n-1$ and coefficients $-1, 0, 1$. Therefore, by Lemma \ref{lemm:garsia} we have
\begin{equation*}
\Delta_n >  c_{\lambda} n^{-d} M_{\lambda} ^ {-n}.
\end{equation*}
The result follows.
\end{proof}

Now we are ready to prove Theorem \ref{theo:main_bernoulli}.

\begin{proof}[Proof of Theorem \ref{theo:main_bernoulli}]
To prove this simply note that letting $F$ be the iterated function system generating the Bernoulli convolution. We have by Corollary \ref{coro:garsia}
\begin{equation*}
M_F \leq M_{\lambda}
\end{equation*}
and by the requirement that $\lambda$ is never root of a non-zero polynomial with coefficients $-1$, $0$, $1$ we have
\begin{equation*}
h_F = \log 2 .
\end{equation*}
To see this note that $h_{F, k}$ is defined to be the entropy of
\begin{equation}
\sum_{i=1}^{k}X_i \lambda^{i-1} \label{eq:spell_out_entropy}
\end{equation}
where each of the $X_i$ are i.i.d. with probability $\frac{1}{2}$ of being each of $\pm 1$. The requirement that $\lambda$ is never root of a non-zero polynomial with coefficients $-1$, $0$, $1$ ensures that each possible choice of the values for the $X_i$ gives a different value for \eqref{eq:spell_out_entropy}. Hence $h_{F, k} = k \log 2$ and so $h_F = \log 2$. We are now done by applying Theorem \ref{theo:main_result}.
\end{proof}

\begin{rema} \label{rema:mahler_bigger_than_2}
We now explain how the requirement that $\lambda$ is not the root of a polynomial with coefficients $0, \pm 1$ forces $M_{\lambda} \geq 2$. This is because $\sum_{i=0}^{n-1} \pm \lambda^i$ is supported on $2^ n$ points each of which are contained in the interval $[-(1-\lambda)^{-1}, (1-\lambda)^{-1}]$. Hence there must be two points in the support with distance at most $2^{-n + o(n)}$. By Lemma \ref{lemm:garsia} it follows that $M_{\lambda} \geq 2$.
\end{rema}
\section{Examples of absolutely continuous self-similar measures} \label{section:examples}
In this section, we give examples of self-similar measures that can be shown to be absolutely continuous using the results of this paper.
\subsection{Examples of absolutely continuous Bernoulli convolutions}
In this subsection, we give explicit values of $\lambda$ for which the Bernoulli convolution with parameter $\lambda$ satisfies the conditions of Theorem \ref{theo:main_bernoulli}.  We do this by a simple computer search. We can ensure that $\lambda$ is not a root of a non-zero polynomial with coefficients $0, \pm 1$ by ensuring that it has a conjugate with absolute value greater than $2$.

The computer search works by checking each integer polynomial with at most a given degree, with all coefficients having at most a given absolute value, with leading coefficient $1$ and with constant term $\pm 1$. The program then finds the roots of the polynomial. If there is one real root with modulus at least $2$ and at least one real root in $(\frac{1}{2}, 1)$, the program then checks that the polynomial is irreducible. If the polynomial is irreducible it then tests each real root in $(\frac{1}{2}, 1)$ to see if it satisfies equation \eqref{eq:main_bernoulli}. In Table \ref{tab:examples} are the results for polynomials of degree at most $11$ and with coefficients of absolute value at most $3$.

The smallest value of $\lambda$ which we were able to find for which the Bernoulli convolution with parameter $\lambda$ can be shown to be absolutely continuous using this method is  $\lambda \approx  0.78207$ with minimal polynomial $X^8 - 2 X^7 - X + 1$.

We were also able to find an infinite family of $\lambda$ for which the results of this paper show that the Bernoulli convolution with parameter $\lambda$ is absolutely continuous. This family is found using the following lemma.

\begin{table}[htb] 
\footnotesize

\begin{tabular}{ | c | c| c| } 

\hline
Minimal polynomial & Mahler measure & $\lambda$\\
\hline
$X^{7} - X^{6} - 2X^{5} - X^{2} + X + 1$ & 2.01043 & 0.87916\\
$X^{7} + 2X^{6} - X - 1$ & 2.01516 & 0.93286\\
$X^{8} - 2X^{7} - X + 1$ & 2.00766 & 0.78207\\
$X^{8} - X^{7} - 2X^{6} - X^{3} + X + 1$ & 2.02530 & 0.90705\\
$X^{8} + 2X^{7} - 1$ & 2.00761 & 0.86058\\
$X^{8} + 2X^{7} + X^{6} + 2X^{5} - X^{2} - X - 1$ & 2.01799 & 0.87735\\
$X^{9} - 2X^{8} - X^{2} + 1$ & 2.01137 & 0.84164\\
$X^{9} - 2X^{8} - X + 1$ & 2.00386 & 0.79953\\
$X^{9} + 2X^{8} - X - 1$ & 2.00386 & 0.94956\\
$X^{9} + 2X^{8} + X^{7} + 2X^{6} - X^{3} - 2X^{2} - X - 1$ & 2.04146 & 0.96868\\
$X^{10} - 2X^{9} - X^{2} + 1$ & 2.00575 & 0.85258\\
$X^{10} - 2X^{9} - X + 1$ & 2.00194 & 0.81397\\
$X^{10} - 2X^{9} + X^{8} - 2X^{7} - X^{5} + X^{4} - X^{3} + 2X^{2} - X + 1$ & 2.02576 & 0.91295\\
$X^{10} - X^{9} - 2X^{8} - X^{7} + X^{6} + 2X^{5} - X^{3} - X^{2} + 1$ & 2.01560 & 0.85694\\
$X^{10} - X^{9} - 2X^{8} - X^{5} + X^{4} + X^{3} - X^{2} + 1$ & 2.01418 & 0.91102\\
$X^{10} - X^{9} - X^{8} - 2X^{7} - X^{5} + X^{4} + X^{2} + 1$ & 2.01224 & 0.93921\\
$X^{10} - X^{9} - X^{8} - X^{7} - 2X^{6} - X^{5} + X^{3} + X^{2} + X + 1$ & 2.01757 & 0.95395\\
$X^{10} - 2X^{8} - 3X^{7} - 2X^{6} - X^{5} + X^{3} + 2X^{2} + 2X + 1$ & 2.00826 & 0.96846\\
$X^{10} + X^{9} - 2X^{8} + X^{7} + X^{6} - X^{5} + X^{4} - X^{3} + X - 1$ & 2.01606 & 0.87581\\
$X^{10} + 2X^{9} - X^{6} - X^{5} + X^{4} - 1$ & 2.03336 & 0.93639\\
$X^{10} + 2X^{9} - X^{4} - 1$ & 2.03066 & 0.94693\\
$X^{10} + 2X^{9} - 1$ & 2.00194 & 0.88881\\
$X^{10} + 3X^{9} + 3X^{8} + 3X^{7} + 2X^{6} - 2X^{4} - 3X^{3} - 3X^{2} - 2X - 1$ & 2.04716 & 0.98447\\
$X^{11} - 2X^{10} - X^{2} + 1$ & 2.00290 & 0.86182\\
$X^{11} - 2X^{10} - X + 1$ & 2.00097 & 0.82615\\
$X^{11} - X^{10} - 2X^{9} - X^{8} + X^{7} + 2X^{6} + X^{5} - X^{4} - 2X^{3} - X^{2} + X + 1$ & 2.00073 & 0.87666\\
$X^{11} - X^{10} - X^{9} - 2X^{8} - X^{4} + X^{2} + X + 1$ & 2.00498 & 0.95290\\
$X^{11} - X^{10} - X^{9} - X^{8} - X^{7} - 2X^{6} - X^{5} + X + 1$ & 2.01424 & 0.83556\\
$X^{11} + X^{10} - 2X^{9} + X^{8} + X^{7} - 2X^{6} + X^{5} + X^{4} - 2X^{3} + X^{2} + X - 1$ & 2.00073 & 0.83139\\
$X^{11} + X^{10} - X^{9} + 2X^{8} + X^{4} - X^{2} + X - 1$ & 2.00498 & 0.80600\\
$X^{11} + 2X^{10} - X - 1$ & 2.00097 & 0.95961\\
$X^{11} + 2X^{10} + X^{2} - 1$ & 2.00290 & 0.81038\\
$X^{11} + 2X^{10} + X^{9} + 2X^{8} - X^{5} - X^{4} - X^{3} - X^{2} - 1$ & 2.03885 & 0.97258\\

\hline
\end{tabular}
\caption{Examples of parameters of Bernoulli convolutions for which Theorem \ref{theo:main_bernoulli} applies} \label{tab:examples}
\end{table}

\begin{lemm} \label{lemm:infinite_family}
Suppose that $n \geq 5$ is an integer and let
\begin{equation*}
p(X) = X^n - 2 X^{n-1} -X +1.
\end{equation*}
Then $p$ has exactly one root in the interval $(\left( \frac{1}{2} \right) ^ {\frac{2}{\sqrt{n-1}}}, 1)$, exactly one root in the interval $(2, 2+2^{2-n})$ and all of the remaining roots are contained in the interior of  the unit disk. Furthermore $p$ is irreducible.
\end{lemm}

Before proving this we need the following result.

\begin{theo}[Rouch\'e's theorem]
Let $f$ and $g$ be holomorphic functions $\mathbb{C} \to \mathbb{C}$ and let $r>0$. Suppose that for all $z\in \mathbb{C}$ such that $|z| = r$ we have
\begin{equation*}
|g(z)| < |f(z)|.
\end{equation*}
Then $f$ and $f+g$ have the same number of zeros with modulus less than $r$.
\end{theo}

\begin{proof}
This is well known. For a proof see for example \cite[Corollary 5.17]{MARSHALL_2019}.
\end{proof}

We are now ready to prove Lemma \ref{lemm:infinite_family}.

\begin{proof}
First we use Rouch\'e's Theorem to prove that all but one of the roots of $p$ is contained in the unit disk. We apply Rouch\'e's Theorem in the form stated above with $f(z) = - 2 z^{n-1} +1$, $g(z) = z^n-z$ and $r =   \left( \frac{1}{2} \right) ^ {\frac{1}{2n-2}}$. A trivial computation which is left to the reader shows that when $|z|=r$ we have $|f(z)| > |g(z)|$. Hence all but one of the roots of $p$ are contained in the ball of radius $\left( \frac{1}{2} \right) ^ {\frac{1}{2n-2}}$.

The other roots can be found by using the intermediate value theorem. Trivial computations show that $p(2) < 0$ and $p(2+2^{2-n}) > 0$. We can also easily compute that $p(1) < 0$ and it is easy to show that $p\left( \left( \frac{1}{2} \right) ^ {\frac{2}{\sqrt{n-1}}} \right) > 0$ whenever $n \geq 5$. Hence there is a root in the interval $(\left( \frac{1}{2} \right) ^ {\frac{2}{\sqrt{n-1}}}, 1)$. In-fact it must be in the interval $(\left( \frac{1}{2} \right) ^ {\frac{2}{\sqrt{n-1}}}, \left( \frac{1}{2} \right) ^ {\frac{1}{2n-2}})$.

The fact that $p$ has only one root in the interval $(\left( \frac{1}{2} \right) ^ {\frac{2}{\sqrt{n-1}}}, 1)$ follows from the fact that it has only one root in the interval $(0, 1)$. Indeed $p'(0) < 0$ and for $x \in (0, 1)$ we have $p''(x) < 0$ hence $p$ is strictly decreasing on $(0, 1)$ and so has at most one root contained in $(0, 1)$.

The fact that $p$ is irreducible follows from the fact that it is a monic integer polynomial with non-zero constant coefficient and all but one of its zero contained in the interior of the unit disk. If $p$ were not irreducible, then one of its factors would need to have all of its roots contained in the interior of the unit disk. This would mean that the product of the roots of this factor would not be an integer, which is a contradiction.
\end{proof}

We now simply let $\lambda_n$ be the root of  $X^n - 2X^{n-1} - X + 1$ contained in the interval $\left( \left( \frac{1}{2} \right) ^{\frac{2}{\sqrt{n-1}}}, 1 \right)$. To show that the Bernoulli convolution with parameter $\lambda_n$ is absolutely continuous using Theorem \ref{theo:main_bernoulli}, it suffices to show that
\begin{equation*}
(\log  (2 + 2^{2-n}) - \log 2) (\log  (2 + 2^{2-n}) )^2 <  \frac{1}{27} \left( \log (2 ) - \log 2 ^{\frac{2}{\sqrt{n-1}}}\right)^3 2^{-\frac{4}{\sqrt{n-1}}}.
\end{equation*}
The left hand side is decreasing in $n$ and the right hand side is increasing in $n$ and for $n=12$ the left hand side is less than the right hand side so for $n \geq 12$ we know that $\mu_{\lambda_n}$ is absolutely continuous. In Table \ref{tab:examples} we show by computing $\lambda_n$ and $M_{\lambda_n}$ for $n = 8, 9, 10$ and $11$ that in fact $\mu_{\lambda_n}$ is absolutely continuous for $n \geq 8$.
\begin{rema}
It is worth noting that we have $\lambda_n \to 1$ and $M_{\lambda_n} \to 2$ so all but finitely many of these Bernoulli convolutions can be shown to be absolutely continuous by the results of  \cite{VARJU_2019}. Using the results of \cite{VARJU_2019} does however require a significantly higher value of $n$ to work. Indeed it requires $n \geq 10^{65}$.
\end{rema}

\subsection{Other examples in dimension one}
In this subsection we briefly mention some other examples of iterated function systems in dimension one that can be shown to be absolutely continuous by these methods.
\begin{prop}
Let $q$ be a prime number and for $i = 1, \dots, q-1$ let $S_i : x \mapsto \frac{q-1}{q} x + i$. Let $F$ be the iterated function system on $\R^1$ given by
\begin{equation*}
 F =\left ((S_i)_{i=1}^q, \left(\frac{1}{q-1}, \dots, \frac{1}{q-1}\right)\right).
\end{equation*}
Then we have $M_F \leq \log q$, $h_F = \log (q-1)$ and $\lambda = \frac{q-1}{q}$. Furthermore, if $q \geq 17$ then $\mu_F$ is absolutely continuous.
\end{prop}
\begin{proof}
We note that any point in the $k$- step iteration of $F$ must be of the form $u = \sum_{i=0}^{k-1} x_i \left( \frac{q-1}{q} \right) ^i$ with $x_i \in \{1, \dots, q-1\}$. Suppose $u = \sum_{i=0}^{k-1} x_i \left( \frac{q-1}{q} \right) ^i$ and  $v = \sum_{i=0}^{k-1} y_i \left( \frac{q-1}{q} \right) ^i$ are two such points. We note that $q^{k-1} u, q^{k-1} v \in \mathbb{Z}$. Therefore, if $u \neq v$ then $|u-v| \geq q^{-(k-1)}$. This gives $M_F \leq \log q$.

We can also note if $u = v$, then looking at $q^{k-1} u$ and $q^{k-1} v $ mod $q^i$ for $i = 1, \dots, k$ we see that we must have $(x_1, x_2, \dots, x_k) = (y_1, y_2, \dots, y_k) $. Therefore, $F$ has no exact overlaps and consequently $h_F = \log (q-1)$. 

We also note that $\lambda = \frac{q-1}{q}$ follows immediately from the definition of $F$.

To show that $\mu_F$ is absolutely continuous using Theorem \ref{theo:main_result} it is sufficient to check that
\begin{equation*}
(\log q - \log (q-1)) (\log q)^2 < \frac{1}{27}\left(\log q - \log \left(\frac{q}{q-1}\right) \right)^3  \left( \frac{q-1}{q} \right)^2.
\end{equation*}
This is the same as showing that
\begin{equation}
\left( \log \left( 1 + \frac{1}{q-1} \right) \right) < \frac{1}{27} \left( \frac{\log (q-1)}{\log q} \right)^2 (\log (q-1)) \left( \frac{q-1}{q} \right)^2. \label{eq:qinequality}
\end{equation}
The left had side of \eqref{eq:qinequality} is decreasing in $q$ and the right hand side is increasing in $q$. The inequality is satisfied for $q = 17$ and so is satisfied for $q \geq 17$.
\end{proof}

\subsection{Examples in dimension two}

In this section we describe some examples of self-similar measures on $\R^2$ which can be shown to be absolutely continuous using the methods of this paper and which cannot be expressed as the product of self-similar measures on $\R$. This is done by identifying $\R^2$ with $\mathbb{C}$.
\begin{prop} \label{prop:2dexamples}
Let $p$ be a prime number such that $p \equiv 3 \, (\text{mod } 4)$. Let $I_p$ denote the ideal $(p)$ in the ring $\mathbb{Z}[i]$. Note that  this is a prime ideal. Let $a_1, \dots, a_{m}$ be in different cosets of $I_p$. Choose some $\alpha$ of the form $\alpha = \frac{a}{p}$ with $a \in \mathbb{Z}[i] \backslash I_p$ and $|\alpha| < 1$. Let $\lambda = |\alpha|$ and let $U : \R^2 \to \R^2$ be a rotation around the origin by $\arg \alpha$.
For $i = 1, \dots, m$ let
\begin{align*}
S_i : \R^2 & \to \R^2 \\
x & \mapsto \lambda U x+a_i
\end{align*}
and let $F$ be the iterated function system on $\R^2$ given by $F = \left( \left( S_i \right)_{i=1}^{m}, \left( \frac{1}{m}, \dots, \frac{1}{m} \right) \right)$. Then we have $M_F \leq \log p$ and $h_F =  \log m$. 
\end{prop}
\begin{proof}
Note that if we identify $\R^2$ with $\mathbb{C}$ then we have
\begin{equation*}
 S_i:z \mapsto \alpha z + a_i.
\end{equation*}
To see that $M_F \leq \log p$ let $x = \sum_{i=0}^{k-1} x_i \alpha^i $ and $y = \sum_{i=0}^{k-1} y_i \alpha^i $ be two points in the $k$-step support of $F$. Note that $p^{k-1} (x-y) \in \mathbb{Z}[i]$ and so if $x \neq y$ then $|x-y| \geq p^{-k+1}$.
To prove $h_F = \log m$ it suffices to show that $F$ has no exact overlaps. For this it suffices to show that if $x_1, \dots, x_k, y_1, \dots, y_k \in \{a_0, \dots, a_{m}\}$ and
\begin{equation}
\sum_{i=0}^{k} x_i \alpha^i = \sum_{i=0}^{k} y_i \alpha^i \label{eq:primeidealthing}
\end{equation} 
then $x_i = y_i$ for $i=1, \dots, k$. We prove this by induction on $i$. For $i=k$ simply multiply both sides of \eqref{eq:primeidealthing} by $p^k$ and then work modulo the ideal $I_p$. Doing this we deduce that $x_k$ and $y_k$ must be in the same coset of $I_p$ which in particular means that they must be equal. The inductive step follows by the same argument.
\end{proof}
		
Note that the above proposition combined with Theorem \ref{theo:main_result} makes it very easy to give numerous examples of absolutely continuous iterated function systems in $\R^2$ which are not products of absolutely continuous iterated function systems in $\R^1$. Some possible examples are given in the following corollary.
\begin{coro}
Let $p$ be a prime number such that $p \equiv 3 \, (\text{mod } 4)$. Let $I_p$ denote the ideal $(p)$ in the ring $\mathbb{Z}[i]$.  Let $a_1, \dots, a_{m}$ be in different cosets of $I_p$. Choose some $\alpha$ of the form $\alpha = \frac{p-1 + i}{p}$. Let $\lambda = |\alpha|$ and let $U : \R^2 \to \R^2$ be a rotation around the origin by $\arg \alpha$.
For $i = 1, \dots, m$ let $S_i :\R^2 \to \R^2, \, x \mapsto \lambda U x+a_i$ and let $F$ be the iterated function system on $\R^2$ given by $F = \left(  \left( S_i \right)_{i=1}^{m}, \left( \frac{1}{m}, \dots, \frac{1}{m} \right) \right)$. Suppose that
\begin{equation*}
(2 \log p - \log m) (\log p)^2 < \frac{1}{27} \left(\log p - \log \frac{p}{p-1} \right)^3 \left( \frac{p-1}{p} \right)^2
\end{equation*}
then the self-similar measure $\mu_F$ is absolutely continuous.
\end{coro}
\begin{proof}
This follows immediately from Theorem \ref{theo:main_result} and Proposition \ref{prop:2dexamples}. Note that in the notation of Theorem \ref{theo:main_result} we have $\lambda \geq \frac{p-1}{p}$.
\end{proof}
\begin{rema}
It is worth noting that the case $m=p^2$ follows from the methods of Garsia \cite{GARSIA_1962}, so in this case the result of this paper can again be seen as a strengthening of the results of \cite{GARSIA_1962}. It is also worth noting that in the case $m=p^2-1$ the conditions of this corollary are satisfies for all $p$ with $p \equiv 3 \, (\text{mod } 4)$ and $p \geq 7$.
\end{rema}

\section{Acknowledgements}

First of all I would like to thank my supervisor Peter  Varj\'u for his help and detailed comments in preparing this paper. I would also like to thank Ioannis Kontoyiannis and Lampros Gavalakis for our discussions on entropy which have helped lead to a more elegant proof of Proposition \ref{prop:entropy_to_detail}. I am grateful to the referees and editors for their careful reading of this paper and for their helpful suggestions which greatly improved its presentation.

\bibliography{refs}{}

\begin{bibdiv}
\begin{biblist}

\bib{breuillard_varju_2020}{article}{
      author={Breuillard, Emmanuel},
      author={Varj{\'u}, P{\'e}ter~P},
       title={Entropy of bernoulli convolutions and uniform exponential growth
  for linear groups},
        date={2020},
     journal={Journal d'Analyse Math{\'e}matique},
      volume={140},
      number={2},
       pages={443\ndash 481},
}

\bib{ERDOS_1939}{article}{
      author={Erd\H{o}s, Paul},
       title={On a family of symmetric {B}ernoulli convolutions},
        date={1939},
     journal={Amer. J. Math.},
      volume={61},
       pages={974\ndash 976},
}

\bib{ERDOS_1940}{article}{
      author={Erd\H{o}s, Paul},
       title={On the smoothness properties of a family of {B}ernoulli
  convolutions},
        date={1940},
        ISSN={0002-9327},
     journal={Amer. J. Math.},
      volume={62},
       pages={180\ndash 186},
         url={https://doi.org/10.2307/2371446},
      review={\MR{858}},
}

\bib{GARSIA_1962}{article}{
      author={Garsia, Adriano~M.},
       title={Arithmetic properties of {B}ernoulli convolutions},
        date={1962},
     journal={Trans. Amer. Math. Soc.},
      volume={102},
       pages={409\ndash 432},
}

\bib{HOCHMAN_2014}{article}{
      author={Hochman, Michael},
       title={On self-similar sets with overlaps and inverse theorems for
  entropy},
        date={2014},
     journal={Ann. of Math. (2)},
      volume={180},
      number={2},
       pages={773\ndash 822},
}

\bib{MR3966837}{inproceedings}{
      author={Hochman, Michael},
       title={Dimension theory of self-similar sets and measures},
        date={2018},
   booktitle={Proceedings of the {I}nternational {C}ongress of
  {M}athematicians---{R}io de {J}aneiro 2018. {V}ol. {III}. {I}nvited
  lectures},
   publisher={World Sci. Publ., Hackensack, NJ},
       pages={1949\ndash 1972},
      review={\MR{3966837}},
}

\bib{HUTCHINSON_1981}{article}{
      author={Hutchinson, John~E.},
       title={Fractals and self-similarity},
        date={1981},
     journal={Indiana Univ. Math. J.},
      volume={30},
      number={5},
       pages={713\ndash 747},
}

\bib{JESSEN_WINTNER_1935}{article}{
      author={Jessen, Borge},
      author={Wintner, Aurel},
       title={Distribution functions and the {R}iemann zeta function},
        date={1935},
     journal={Trans. Amer. Math. Soc.},
      volume={38},
      number={1},
       pages={48\ndash 88},
}

\bib{Johnson_2004}{book}{
      author={Johnson, Oliver~T.},
       title={Information theory and the central limit theorem},
   publisher={World Scientific Publishing Company},
        date={2004},
}

\bib{MR1507093}{article}{
      author={Kershner, Richard},
      author={Wintner, Aurel},
       title={On {S}ymmetric {B}ernoulli {C}onvolutions},
        date={1935},
        ISSN={0002-9327},
     journal={Amer. J. Math.},
      volume={57},
      number={3},
       pages={541\ndash 548},
         url={https://doi.org/10.2307/2371185},
      review={\MR{1507093}},
}

\bib{KITTLE_2022}{article}{
      author={Kittle, Samuel},
       title={Absolutely continuous furstenberg measures with exponetial
  separation},
     journal={work in progress},
}

\bib{Madiman_Kontoyiannis_2014}{article}{
      author={Kontoyiannis, Ioannis},
      author={Madiman, Mokshay},
       title={Sumset and inverse sumset inequalities for differential entropy
  and mutual information},
        date={2014},
        ISSN={0018-9448},
     journal={IEEE Trans. Inform. Theory},
      volume={60},
      number={8},
       pages={4503\ndash 4514},
         url={https://doi.org/10.1109/TIT.2014.2322861},
      review={\MR{3245338}},
}

\bib{MARSHALL_2019}{book}{
      author={Marshall, Donald~E.},
       title={Complex analysis},
      series={Cambridge Mathematical Textbooks},
   publisher={Cambridge University Press, Cambridge},
        date={2019},
        ISBN={978-1-107-13482-9},
      review={\MR{4321146}},
}

\bib{PERES_SCHLAG_WILHEM_2000}{incollection}{
      author={Peres, Yuval},
      author={Schlag, Wilhelm},
      author={Solomyak, Boris},
       title={Sixty years of {B}ernoulli convolutions},
        date={2000},
   booktitle={Fractal geometry and stochastics ii},
      series={Progr. Probab.},
      volume={46},
   publisher={Birkh\"{a}user, Basel},
       pages={39\ndash 65},
}

\bib{SHMERKIN_2014}{article}{
      author={Shmerkin, Pablo},
       title={On the exceptional set for absolute continuity of {B}ernoulli
  convolutions},
        date={2014},
     journal={Geom. Funct. Anal.},
      volume={24},
      number={3},
       pages={946\ndash 958},
}

\bib{SHMERKIN_2019}{article}{
      author={Shmerkin, Pablo},
       title={On {F}urstenberg's intersection conjecture, self-similar
  measures, and the {$L^q$} norms of convolutions},
        date={2019},
     journal={Ann. of Math. (2)},
      volume={189},
      number={2},
       pages={319\ndash 391},
}

\bib{SOLOMYAK_1995}{article}{
      author={Solomyak, Boris},
       title={On the random series {$\sum\pm\lambda^n$} (an {E}rd{\H{o}}s
  problem)},
        date={1995},
        ISSN={0003-486X},
     journal={Ann. of Math. (2)},
      volume={142},
      number={3},
       pages={611\ndash 625},
         url={https://doi.org/10.2307/2118556},
      review={\MR{1356783}},
}

\bib{SOLOMYAK_2004}{incollection}{
      author={Solomyak, Boris},
       title={Notes on {B}ernoulli convolutions},
        date={2004},
   booktitle={Fractal geometry and applications: a jubilee of {B}eno\^{\i}t
  {M}andelbrot. {P}art 1},
      series={Proc. Sympos. Pure Math.},
      volume={72},
   publisher={Amer. Math. Soc., Providence, RI},
       pages={207\ndash 230},
      review={\MR{2112107}},
}

\bib{VARJU_2021}{incollection}{
      author={Varj\'{u}, P\'{e}ter~P.},
       title={Recent progress on {B}ernoulli convolutions},
        date={2018},
   booktitle={European {C}ongress of {M}athematics},
   publisher={Eur. Math. Soc., Z\"{u}rich},
       pages={847\ndash 867},
      review={\MR{3890454}},
}

\bib{VARJU_2019}{article}{
      author={Varj\'{u}, P\'{e}ter~P.},
       title={Absolute continuity of {B}ernoulli convolutions for algebraic
  parameters},
        date={2019},
        ISSN={0894-0347},
     journal={J. Amer. Math. Soc.},
      volume={32},
      number={2},
       pages={351\ndash 397},
         url={https://doi.org/10.1090/jams/916},
      review={\MR{3904156}},
}

\end{biblist}
\end{bibdiv}
\bibliographystyle{smfplain}

\end{document}